\documentclass[MSNbibl,number,citesort,rotating,dvips]{arxbj}
\usepackage{dcolumn,upgreek,mathbbol}
\usepackage{graphicx}

\newcommand{\bbSigma}{\bold{\mathbb{\Sigma}}}

\aid{0}
\volume{19}
\issue{1}
\pubyear{2013}
\firstpage{252}
\lastpage{274}
\doi{10.3150/11-BEJ386}

\makeatletter
\define@key{bid}{pmid}{}

\def\P{\mathrm{P}}
\def\boldmu{\bolds{\mu}}
\newcommand{\sgn}{\operatorname{sgn}}
\newcommand{\diag}{\operatorname{diag}}
\newtheorem{theorem}{Theorem}
\newtheorem{lemma}{Lemma}
\newremark{remark}{Remark}
\newcolumntype{d}[1]{D{.}{.}{#1}}
\newcommand{\T}{\mathrm{T}}
\makeatother

\begin{document}
\begin{frontmatter}

\title{Simultaneous variable selection and estimation in semiparametric
modeling of longitudinal/clustered data}
\runtitle{Penalized semiparametric regression for clustered data}

\begin{aug}
\author[a]{\fnms{Shujie} \snm{Ma}\thanksref{a}\ead[label=e1]{shujie.ma@ucr.edu}},
\author[b]{\fnms{Qiongxia} \snm{Song}\thanksref{b}\ead[label=e2]{song@utdallas.edu}}
\and
\author[c]{\fnms{Li} \snm{Wang}\corref{}\thanksref{c}\ead[label=e3]{lilywang@uga.edu}}
\runauthor{S. Ma, Q. Qiongxia and L. Wang} 
\address[a]{Department of Statistics, UC-Riverside, Riverside, CA
92521, USA.
\printead{e1}}
\address[b]{Department of Mathematical Sciences, University of Texas at
Dallas, Richardson, TX 75080, USA.
\printead{e2}}
\address[c]{Department of Statistics, The University of Georgia,
Athens, GA 30602,
USA.\\
\printead{e3}}
\end{aug}

\received{\smonth{8} \syear{2010}}
\revised{\smonth{6} \syear{2011}}

%
\begin{abstract}
We consider the problem of simultaneous variable selection and
estimation in additive, partially linear models for
longitudinal/clustered data. We propose an estimation procedure via
polynomial splines to estimate the nonparametric components and apply
proper penalty functions to achieve sparsity in the linear part. Under
reasonable conditions, we obtain the asymptotic normality of the
estimators for the linear components and the consistency of the
estimators for the nonparametric components. We further demonstrate
that, with proper choice of the regularization parameter, the penalized
estimators of the non-zero coefficients achieve the asymptotic oracle
property. The finite sample behavior of the penalized estimators is
evaluated with simulation studies and illustrated by a longitudinal CD4
cell count data set.
\end{abstract}

%

\end{frontmatter}

\section{Introduction}\label{SEC:Introduction}

In the past two decades, there has been a considerable amount of
research to study additive, partially linear models (APLM); see Opsomer and Ruppert \cite
{OR99}, H\"{a}rdle, Liang and Gao~\cite{HLG00}, Li~\cite{L00}, Fan and Li~\cite{FL03}, Liang \textit{et al.}~\cite{LTRAH08}, Liu, Wang and Liang \cite
{LWL11}, Ma and Yang~\cite{MY11}, among others. APLMs meet three fundamental
aspects (Stone~\cite{S85}) of statistical models: flexibility,
dimensionality and interpretability. In this paper, we consider the
APLMs for clustered and longitudinal data.

Let $\{(Y_{ij},\mathbf{X}_{ij},\mathbf{Z}_{ij}),1\leq i\leq
n,1\leq j\leq m_{i}\}$ be the $j$th observation for the $i$th
subject or cluster, where $Y_{ij}$ is the response variable, $%
\mathbf{X}_{ij}=( 1,X_{ij1},\ldots,X_{ij(d_{1}-1)}) ^{\T}$ is a $%
d_{1}$-vector of covariates, and $\mathbf{Z}_{ij}=(
Z_{ij1},\ldots,Z_{ijd_{2}}) ^{\T}$ is a $d_{2}$-vector of
covariates. An APLM for this kind of data is given by
\begin{equation}
Y_{ij}=\mu_{ij}+\varepsilon_{ij}=\mathbf{X}_{ij}^{\T}\bolds{\beta}%
+\sum_{l=1}^{d_{2}}\eta_{l}(Z_{ijl})+\varepsilon_{ij},\qquad j=1,\ldots
, m_{i},i=1,\ldots,n, \label{model}
\end{equation}
where $\bolds{\beta}$ is a $d_{1}$-dimensional regression parameter,
and $%
\eta_{l}$, $l=1,\ldots,d_{2}$, are unknown but smooth functions. We
assume $%
\underline{\bolds{\varepsilon}}_{i}=( \varepsilon_{i1},\ldots
,\varepsilon_{im_{i}}) ^{\T}\sim N( \mathbf{0},\bolds{\Sigma}%
_{i}) $. For identifiability, both the parametric and nonparametric components
must be centered, that is, $E\eta_{l}(Z_{ijl})\equiv0$,
$l=1,\ldots,d_{2}$, $%
EX_{ijk}=0,$ $k=1,\ldots,d_{1}$. When $d_{2}=1$, model (\ref{model}) is
simplified to be the partially linear model (PLM) in Lin and Carroll~\cite{LC01}. Model
(\ref%
{model}) retains the merits of additive models, while it is more flexible
than purely additive models by allowing a subset of the covariates to
be discrete and/or unbounded. When $m_{i}$s and $\bolds{%
\Sigma}_{i}$s are the same for all individuals, Carroll \textit{et al.} \cite
{CMMY09} considered
the efficient estimation of $\bolds{\beta}$ in model (\ref{model}) using
local linear smooth backfitting. In this paper we consider a more general
scenario that both $m_{i}$ and $\bolds{\Sigma}_{i}$ may vary across
subjects or experimental units to allow irregular measurements for
individuals. Our goal is to simultaneously select significant variables
and efficiently estimate the unknown components for model (\ref
{model}). This is challenging due to
the issue of ``curse of dimensionality'' and the
additional complexity of the correlation structures (Wang~\cite{W03})
introduced by repeated measurements.

To alleviate the effect of the ``curse of
dimensionality,'' more parsimonious models become desirable in
practice; see Fan~\cite{F03}, Hall, M\"{u}ller and Wang~\cite{HMW06} and Wang \textit{et al.}~\cite{WXZC}. Variable
selection is fundamental to
high-dimensional statistical modeling. In the absence of prior
knowledge, a
large number of variables may be included at the initial stage of modeling
in order to reduce possible model bias. This may lead to a complicated model
including many insignificant variables, resulting in less predictive powers
and difficulty in interpretation. There is an extensive literature on
variable selection via various approaches, for example, the classical
information criteria such as the Akaike information criterion (AIC) and
Bayesian information criterion (BIC) in Yang~\cite{Y08}, the least absolute
shrinkage and selection operator (LASSO) proposed in Tibshirani~\cite%
{T96,T97}, the non-negative garrote in Yuan and Liu~\cite{YL07}, the difference
convex algorithm in Wu and Liu~\cite{WL09}, the combination of $L_{0}$ and $L_{1}$
penalties in Liu and Wu~\cite{LW07}, and the nonparametric independence screening
procedure in Fan, Feng and Song~\cite{FFS11}.

Many traditional variable selection procedures in use, including stepwise
selection, AIC or BIC, can be expensive in computation and ignore stochastic
errors inherited in the variable selection process. Penalized least
squares approaches have gained popularity in recent years to
automatically and simultaneously select significant variables; for
example, Antoniadis
\cite{A97} proposed the hard thresholding penalty which enables best
subset selection and stepwise deletion in certain cases. The LASSO
(Tibshirani~\cite{T96,T97}) is one of the most popular shrinkage estimators,
but it has some deficiencies (Meinshausen and B\"{u}hlmann~\cite{MB06}). Fan and Li~\cite{FL01}
proposed the smoothly clipped absolute deviation penalty (SCAD), which
achieves an ``oracle'' property in the sense that it performs as well
as if the subset of significant variables were known in advance. The
SCAD-penalized selection procedures were illustrated in Fan and Li \cite
{FL01} for parametric models; Cai \textit{et al.}~\cite{CFLZ05} and Fan and Li~\cite{FL02}
for survival models; Li and Liang~\cite{LL08} for generalized varying-coefficient
models; Liang and Li~\cite{LL09}\vadjust{\goodbreak} and Ma and Li~\cite{ML10} for measurement error models; Xue~\cite
{X09} for pure additive models; and Xue, Qu and Zhou~\cite{XQZ10} for generalized
additive models with correlated data.

We propose a model selection method for APLMs with repeated measures by
penalizing appropriate estimating functions. We approximate
nonparametric components by spline functions and obtain asymptotic
normality for the coefficient estimators via one step least squares.
The proposed approach is computationally expedient and easy to
implement, in contrast to the backfitting approach in Carroll \textit{et al.}
\cite{CMMY09}. Moreover, it avoids the pitfall of the
backfitting algorithms caused by dependence between covariates.
Furthermore, we show that the estimator can correctly select the
nonzero coefficients with probability converging to $1$ and the $\sqrt
{n}$-consistent estimators of the non-zero coefficients can perform as
well as an oracle estimator in the sense of Fan and Li~\cite{FL01} with
a suitable choice of penalty function.

The paper is organized as follows. In Section~\ref{SEC:Methodology}, we
introduce the penalized polynomial spline estimating method. Section
\ref{SEC:Asymptotics} provides the asymptotic properties of the proposed
estimators, including the consistency and oracle property of the parametric
components, as well as the rate of the $L_{2}$-convergence of the
nonparametric components. In Section~\ref{SEC:Implementation}, we
discuss some implementation issues of the proposed procedure.
Simulation studies are
presented in Section~\ref{SEC:Simulation}. Section~\ref{SEC:Application}
illustrates the application using longitudinal CD4 cell-count data. We
conclude with a discussion in Section~\ref{SEC:Discussion}. Technical
proofs are presented in the \hyperref[appendix]{Appendix}.

\section{Penalized spline estimation}
\label{SEC:Methodology}

For simplicity, denote vectors $\underline{\mathbf{Y}}_{i}= (
Y_{i1},\ldots,Y_{im_{i}}) ^{\T} $ and $%
\underline{\bolds{\mu}}_{i}= ( \mu_{i1},\ldots,\mu
_{im_{i}}) ^{\T}$, $1\leq m_{i}\leq M$, $%
1\leq i\leq n$. Similarly, let $\underline{\mathbf{X}}_{i}=\{ (
\mathbf{X}_{i1},\ldots,\mathbf{X}_{im_{i}}) ^{\T}\}
_{m_{i}\times d_{1}}$ and $\underline{\mathbf{Z}}_{i}=\{ ( \mathbf{Z}_{i1},\ldots,\mathbf{Z}_{im_{i}}) ^{\T}\} _{m_{i}\times
d_{2}} $. Assume that $Z_{ijl}$ has the same distribution as $Z_{l}$, which
is distributed on a compact interval $[ a_{l},b_{l}] ,1\leq l\leq
d_{2}$, and, without loss of generality, we take all intervals $[
a_{l},b_{l}] =[ 0,1] ,1\leq l\leq d_{2}$. Let $\eta_{l}(%
\mathbf{Z}_{il})=\{ \eta_{l}(Z_{i1l}),\ldots,\eta
_{l}(Z_{im_{i}l})\} ^{\T}$, for $%
l=1,\ldots,d_{2} $. The mean function in model~(\ref{model}) can be
written in
matrix notation as $\underline{\bolds{\mu}}_{i}=\bolds{\underline
{X}}_{i}%
\bolds{\beta}+\sum_{l=1}^{d_{2}}\eta_{l}(\mathbf{Z}_{il})$,
which is a semiparametric extension of the marginal model in Liang and Zeger \cite
{LZ86} with an identity link.

As in Wang, Carroll and Lin~\cite{WCL05}, we allow $\mathbf{X}$ and $\mathbf{Z}$ to be dependent.
Let $\mathbf{V}_{i}=\mathbf{V}_{i}( \underline{\mathbf{X}}_{i},\underline{\mathbf{Z}}_{i}) $ be the assumed ``working''\ covariance of
$\underline{\mathbf{Y}}_{i}$, where
$\mathbf{V}_{i}=\mathbf{A}_{i}^{1/2}\mathbf{R}_{i}\mathbf{A}_{i}^{1/2}$, $%
\mathbf{A}_{i}$ denotes a $m_{i}\times m_{i}$ diagonal matrix that contains
the marginal variances of $Y_{ij}$, and $\mathbf{R}_{i}$ is an invertible
working correlation matrix. Throughout, we assume that $\mathbf{V}_{i}$
depends on a nuisance finite dimensional parameter vector $\bolds
{\alpha}$.

Following Wang and Yang~\cite{WY07}, we approximate the nonparametric functions $\eta
_{l}$'s by polynomial splines. Let $G_{n}$ be the space of polynomial
splines of degree $q\geq1$. We introduce a sequence of spline knots
\[
t_{-q}=\cdots=t_{-1}=t_{0}=0<t_{1}<\cdots<t_{N}<1=t_{N+1}=\cdots=t_{N+q+1},
\]
where $N\equiv N_{n}$ is the number of interior knots, and $N$
increases when sample size $n$ increases with the
precise order given in Assumption (A5).\vadjust{\goodbreak} Then $G_{n}$ consists of
functions $%
\varpi$ satisfying (i)~$\varpi$ is a polynomial of degree $q$ on each of
the subintervals $I_{s}=[ t_{s},t_{s+1}) $, $s=0,\ldots,N_{n}-1$, $%
I_{Nn}=[ t_{N_{n}},1] $; (ii) for $q\geq1$, $\varpi$ is $(q-1)$
times continuously differentiable on $[0,1]$. In the following, let
$J_{n}=N_{n}+q+1$, and we adopt the normalized B-spline space
$G_{n}^{0}=\{ B_{s,l}\dvt 1\leq l\leq
d_{2},1\leq s\leq J_{n}\} ^{\T}$ in Xue and Yang~\cite{XY06}. Equally spaced knots
are used in this article for simplicity of proof. However, other
regular knot
sequences can also be used with similar asymptotic results.

Suppose that $\eta_{l}$ can be approximated well by a spline function
in $%
G_{n}^{0}$ so that
\begin{equation}
\eta_{l}( z_{l}) \approx\widetilde{\eta}_{l}(
z_{l}) =\sum_{s=1}^{J_{n}}\gamma_{sl}B_{s,l}( z_{l}) .
\label{EQ:spline}
\end{equation}
Let $\bolds{\gamma}=( \gamma_{sl}\dvtx 1\leq s\leq J_{n},1\leq l\leq
d_{2}) ^{\T}$ be the collection of the coefficients in (\ref{EQ:spline}%
), and let
\begin{equation}
\mathbf{B}_{ijl}=[ \{ B_{s,l}( Z_{ijl}) \dvt 1\leq s\leq
J_{n}\} ^{\T}] _{J_{n}\times1},\qquad \mathbf{B}_{ij}=\{
( \mathbf{B}_{ij1}^{\T},\ldots,\mathbf{B}_{ijd_{2}}^{\T}) ^{\T%
}\} _{d_{2}J_{n}\times1}; \label{EQ:B}
\end{equation}
then we have an approximation $\mu_{ij}\approx\mathbf{X}_{ij}^{\T
}\bolds{\beta}+\mathbf{B}_{ij}^{\T}\bolds{\gamma}$. We can also write the
approximation in
matrix notation as $\underline{\bolds{\mu}}_{i}\approx\underline{\mathbf{X}}_{i}\bolds{\beta}+\underline{\mathbf{B}}_{i}\bolds{\gamma}$,
where $%
\underline{\mathbf{B}}_{i}=\{ ( \mathbf{B}_{i1},\ldots,\mathbf{B}%
_{im_{i}}) ^{\T}\} _{m_{i}\times d_{2}J_{n}}$.

Let $\widehat{\bolds{\beta}}=(
\widehat{\beta}_{1},\ldots,\widehat{\beta}_{d_{1}}) ^{\T}$ and $%
\widehat{\bolds{\gamma}}=\{ \widehat{\gamma}%
_{sl}\dvt s=1,\ldots,J_{n},l=1,\ldots,d_{2}\} ^{\T}$ be the minimizer of
\begin{equation}
Q_{n}( \bolds{\beta,\gamma}) =\frac{1}{2}\sum_{i=1}^{n}\{
\underline{\mathbf{Y}}_{i}-( \underline{\mathbf{X}}_{i}\bolds{\beta}+%
\underline{\mathbf{B}}_{i}\bolds{\gamma}) \} ^{\T}\mathbf{V}%
_{i}^{-1}\{ \underline{\mathbf{Y}}_{i}-( \underline{\mathbf{X}}_{i}%
\bolds{\beta}+\underline{\mathbf{B}}_{i}\bolds{\gamma}) \} ,
\label{DEF:Qn}
\end{equation}
which is corresponding to the class of working covariance matrices $\{
\mathbf{V}_{i},1\leq i\leq n\} $, or, equivalently, they solve the
estimating equations
\begin{eqnarray}
\sum_{i=1}^{n}\underline{\mathbf{X}}_{i}^{\T}\mathbf{V}_{i}^{-1}\{
\underline{\mathbf{Y}}_{i}-( \underline{\mathbf{X}}_{i}\bolds{\beta}+%
\underline{\mathbf{B}}_{i}\bolds{\gamma}) \} &=&\mathbf{0},
\label{DEF:Sn} \\
\sum_{i=1}^{n}\underline{\mathbf{B}}_{i}^{\T}\mathbf{V}_{i}^{-1}\{
\underline{\mathbf{Y}}_{i}-( \underline{\mathbf{X}}_{i}\bolds{\beta}+%
\underline{\mathbf{B}}_{i}\bolds{\gamma}) \} &=&\mathbf{0}.
\label{EQ:LSE}
\end{eqnarray}
Solving (\ref{EQ:LSE}) yields
\begin{equation}
\bolds{\gamma}\equiv\bolds{\gamma}( \bolds{\beta}) =\Biggl(
\sum_{i=1}^{n}\underline{\mathbf{B}}_{i}^{\T}\mathbf{V}_{i}^{-1}\mathbf{
\underline{B}}_{i}\Biggr) ^{-1}\sum_{i=1}^{n}\underline{\mathbf{B}}_{i}^{\T}%
\mathbf{V}_{i}^{-1}( \underline{\mathbf{Y}}_{i}-\underline{\mathbf{X}}%
_{i}\bolds{\beta}) . \label{EQ:gamma(beta)}
\end{equation}
Replacing $\bolds{\gamma}$ by $\bolds{\gamma}( \bolds{\beta}%
) $ in (\ref{DEF:Qn}), we define
\begin{eqnarray}\label{DEF:Qbeta}
Q( \bolds{\beta}) \equiv Q_{n}\{ \bolds{\beta,\gamma}%
( \bolds{\beta}) \} &=&\frac{1}{2}\sum_{i=1}^{n}[
\underline{\mathbf{Y}}_{i}-\{ \underline{\mathbf{X}}_{i}\bolds{\beta}+
\underline{\mathbf{B}}_{i}\bolds{\gamma}( \bolds{\beta})
\} ] ^{\T}
\nonumber
\\[-8pt]
\\[-8pt]
\nonumber
&&\phantom{\frac{1}{2}\sum_{i=1}^{n}}{}\times\mathbf{V}_{i}^{-1}[ \underline{\mathbf{Y}}_{i}-\{
\underline{\mathbf{X}}_{i}\bolds{\beta}+\bolds{\underline
{B}}_{i}\bolds{%
\gamma}( \bolds{\beta}) \} ].
\end{eqnarray}
To select the significant parametric components, we add a penalty to
$Q( \bolds{\beta})$. Let $n_{\mathrm{T}}=\sum_{i=1}^{n}m_{i}$, and
define the penalized version of $Q( \bolds{\beta})$ as
\begin{equation}
Q_{\mathcal{P}}( \bolds{\beta}) =Q( \bolds{\beta})
+n_{\mathrm{T}}\mathcal{P}( \bolds{\beta}) , \label{DEF:PLS}
\end{equation}
where $\mathcal{P}( \bolds{\beta})
=\sum_{k=1}^{d_{1}}p_{\lambda_{k}}( \vert\beta_{k}\vert
) $ for a pre-specified penalty function $p_{\lambda}( \vert\beta\vert
) $ with a regularization parameter $
\lambda$. Minimizing $Q_{\mathcal{P}}(\bolds{\beta})$ in (\ref%
{DEF:PLS}) yields a penalized estimator
\begin{equation}
\widehat{\bolds{\beta}}^{\P}=\arg\min Q_{\mathcal{P}}( \bolds{%
\beta}) . \label{DEF:beta-PGEE}
\end{equation}

Various penalty functions can be used for $\mathcal{P}( \bolds{\beta
}) $ in variable selection procedures.
We consider two penalty functions, the hard thresholding penalty
(Antoniadis~\cite{A97}) $p_{\lambda
}( \beta) =\lambda^{2}-( \vert\beta\vert
-\lambda) ^{2}I( \vert\beta\vert<\lambda) $%
and the SCAD penalty (Fan and Li~\cite{FL01}), given by
\[
p_{\lambda}^{\prime}( \beta) =\lambda\biggl\{ I( \beta
\leq\lambda) +\frac{( a\lambda-\beta) _{+}}{(
a-1) \lambda}I( \beta>\lambda) \biggr\} \qquad\mbox{for some }%
a>2\mbox{ and }\beta>0,
\]
where $p_{\lambda}(0)=0$, and $\lambda$ and $a$ are two tuning parameters.
Justifying from a Bayesian statistical point of view, Fan and Li \cite
{FL01} suggested
using $a=3.7$, which will be used in our simulation studies.

The minimization problem in (\ref{DEF:beta-PGEE}) is essentially a
one-step least squares problem, which can be easily solved and
implemented with many existing regression programs. The theorems
established in Section~\ref{SEC:penalizedestimators} demonstrate that
$\widehat{\bolds{\beta}}^{\P}$
performs asymptotically as well as an oracle estimator in terms of selecting
the correct model when the regularization parameter is appropriately chosen.

\section{Asymptotic properties of the estimators}
\label{SEC:Asymptotics}

For positive numbers $a_{n}$ and $b_{n}$, $n\geq1$, let $a_{n}\sim b_{n}$
denote that $\lim_{n\rightarrow\infty}a_{n}/b_{n}=c$, where $c$ is some
non-zero constant. Let $\vert\phi\vert_{L_{2}}\equiv[
\int_{0}^{1}\{ \phi( z) \} ^{2}\,\mathrm{d}z] ^{1/2}$
denote the $L_{2}$ norm of any square integrable function $\phi(
z) $ on $[ 0,1] $.
Denote the space of the $p$th order smooth functions as $C^{(
p) }[ 0,1] =$ $\{ \phi\mid\phi^{(
p) }\in C[ 0,1] \} $.

\subsection{Assumptions}

The assumptions for the asymptotic results are listed below:

\begin{enumerate}[(A1)]
\item[(A1)] The random variables $Z_{ijl}$ are bounded, uniformly in
$1\leq
j\leq m_{i}$, $1\leq i\leq n$, $1 \leq l \leq d_{2}$. The marginal
density $%
f_{l}( z_{l}) $ of $Z_{l}$ has the uniform upper bound $C_{f}$
and lower bound $c_{f}$ on $[ 0,1] $. The joint density $%
f_{ll^{\prime}}( z_{l},z_{l^{\prime}}) $ of $(
Z_{ijl},Z_{ijl^{\prime}}) $ satisfies that $c_{f}\leq f_{ll^{\prime
}}( z_{l},z_{l^{\prime}}) \leq C_{f}$, for all $(
z_{l},z_{l^{\prime}}) \in[ 0,1] ^{2}$, $1\leq l\neq
l^{\prime}\leq d_{2}.$

\item[(A2)] The random variables $X_{ijk}$ are bounded, uniformly in
$1\leq
j\leq m_{i}$, $1\leq i \leq n$, $1\leq k\leq d_{1}$. The eigenvalues of
$%
E\{ \mathbf{X}_{ij}\mathbf{X}_{ij}^{\T}\vert\mathbf{Z}%
_{ij} \} $ are bounded away from $0$ and infinity, uniformly in
$1\leq j\leq m_{i}$, $1\leq i\leq n$.

\item[(A3)] The eigenvalues of the true covariance matrices $\bolds
{\Sigma}_{i}$ are bounded away from $0$ and infinity, uniformly in
$1\leq i\leq n$.

\item[(A4)] The eigenvalues of the working covariance matrices $\mathbf
{V}%
_{i}$ are bounded away from $0$ and infinity, uniformly in $1\leq i\leq n$.
\end{enumerate}

To make $\bolds{\beta}$ estimable at the $\sqrt{n}$ rate, we need a
condition to ensure that $\mathbf{X}$ and $\mathbf{Z}$ not functionally related.
Define $\mathcal{H}=\{ \psi( \mathbf{z})
=\sum_{l=1}^{d_{2}}\psi_{l}( z_{l}) , E\psi
_{l}( z_{l}) =0,\vert\psi_{l}\vert_{L_{2}}<\infty
\} $ the Hilbert space of theoretically centered $L_{2}$ additive
functions on $[ 0,1] ^{d_{2}}$. Let $\psi_{k}^{\ast}$ be the
function $\psi\in\mathcal{H}$ that minimizes
\[
\sum_{i=1}^{n}E\bigl[ \bigl\{ \underline{\mathbf{X}}_{i}^{(
k) }-\psi( \underline{\mathbf{Z}}_{i}) \bigr\} ^{\T}%
\mathbf{V}_{i}^{-1}\bigl\{ \underline{\mathbf{X}}_{i}^{( k)
}-\psi( \underline{\mathbf{Z}}_{i}) \bigr\} \bigr] ,
\]
where
\begin{equation}
\underline{\mathbf{X}}_{i}^{( k) }=( X_{i1k},\ldots,
X_{im_{i}k}) ^{\T},\qquad 1\leq k\leq d_{1}. \label{DEF:X_i(k)}
\end{equation}
Then

\begin{enumerate}[(A5)]
\item[(A5)]
for $1\leq l\leq d_{2}$, $1\leq k\leq d_{1}$, assume that $\eta
_{l}(z_{l})\in C^{( p) }[ 0,1] $, $\psi_{k}^{\ast
}\in C^{( p) }[ 0,1] $ for a given integer $p\geq1$,
and the spline degree satisfies $q+1\geq p$. The number of the spline
basis functions $J_{n}\sim n^{1/(2p)}\log( n)$.
\end{enumerate}

Assumptions (A1)--(A4) are identical with (C1)--(C4) in Huang, Zhang and Zhou~\cite{HZZ06}, while
Assumption (A5) is similar to (C1) and (C4) in Liu, Wang and Liang~\cite{LWL11}.

\subsection{Asymptotic properties for the unpenalized estimators}

According to the equations in (\ref{DEF:Sn}) and (\ref{EQ:LSE}), we have
\begin{equation}
\pmatrix{\widehat{\bolds{\beta}}\cr{\widehat{\bolds{\gamma}}}}=\Biggl(
\sum_{i=1}^{n}\underline{\mathbf{D}}_{i}^{\T}\mathbf{V}_{i}^{-1}\underline{\mathbf{D}}_{i}\Biggr) ^{-1}\Biggl( \sum_{i=1}^{n}\underline{\mathbf{D}}%
_{i}^{\T}\mathbf{V}_{i}^{-1}\underline{\mathbf{Y}}_{i}\Biggr) ,
\label{EQ:betagammahat}
\end{equation}
where $\underline{\mathbf{D}}_{i}=( \underline{\mathbf{X}}_{i},\underline{\mathbf{B}}_{i}) _{m_{i}\times( d_{1}+d_{2}J_{n}) }$.
The centered additive component $\eta_{l}(z_{l})$ is estimated by the
empirically centered estimator
\begin{equation}
\widehat{\eta}_{l}(z_{l})=\sum_{s=1}^{J_{n}}\widehat{\gamma}%
_{sl}B_{s,l}( z_{l}) -n_{\mathrm{T}}^{-1}\sum_{i=1}^{n}%
\sum_{j=1}^{m_{i}}B_{s,l}( Z_{ijl}) . \label{EQ:mhat}
\end{equation}

Next we derive the asymptotic properties of $\widehat{\bolds{\beta}}$
and $%
\widehat{\eta}_{l}$. Let $\mathbb{X}$ and $\mathbb{Z}$ be the collections
of all $X_{ijk}$s and $Z_{ijl}$s, respectively, that is, $\mathbb
{X}_{n_{%
\mathrm{T}}\times d_{1}}=( \underline{\mathbf{X}}_{1}^{\T},\ldots,%
\underline{\mathbf{X}}_{n}^{\T}) ^{\T}$ and $\mathbb{Z}_{n_{\mathrm{T}%
}\times d_{2}}=( \underline{\mathbf{Z}}_{1}^{\T},\ldots,\underline{\mathbf{Z}}_{n}^{\T}) ^{\T}$. Define
\begin{equation}
\widetilde{\mathbf{X}}_{i}^{(k)}=\underline{\mathbf{X}}_{i}^{( k)
}-\psi_{k}^{\ast}( \underline{\mathbf{Z}}_{i}) ,\qquad 1\leq k\leq
d_{1}, \qquad \widetilde{\underline{\mathbf{X}}}_{i}=\bigl( \widetilde{\bolds{%
X}}_{i}^{(1)},\ldots,\widetilde{\mathbf{X}}_{i}^{(d_{1})}\bigr)
_{m_{i}\times d_{1}}, \label{EQ:Xtilda}
\end{equation}
for $1\leq i\leq n$. Denote $\widetilde{\mathbb{X}}=\{ (
\underline{\widetilde{\mathbf{X}}}_{1}^{\T},\ldots,\underline
{\widetilde{%
\mathbf{X}}}_{n}^{\T}) ^{\T}\} _{n_{\mathrm{T}}\times d_{1}}$,
\[
\mathbb{V}^{-1}=\operatorname{diag}( \mathbf{V}_{1}^{-1},\ldots,\mathbf{V}%
_{n}^{-1}) _{n_{\mathrm{T}}\times n_{\mathrm{T}}},\qquad \bbSigma =%
\operatorname{diag}( \bolds{\Sigma}_{1},\ldots,\bolds{\Sigma}_{n})
_{n_{\mathrm{T}}\times n_{\mathrm{T}}}.
\]
Further define
\begin{equation}
\bolds{\Omega}( \mathbb{V},\bbSigma ) =\{
\widetilde{\mathbf{A}}( \mathbb{V}) \} ^{-1}\widetilde{%
\mathbf{B}}( \mathbb{V},\bbSigma ) \{ \widetilde{%
\mathbf{A}}( \mathbb{V}) \} ^{-1} \label{EQ:IV}
\end{equation}
with $\widetilde{\mathbf{A}}( \mathbb{V}) =E( n^{-1}%
\widetilde{\mathbb{X}}^{\T}\mathbb{V}^{-1}\widetilde{\mathbb{X}}) $
and $\widetilde{\mathbf{B}}( \mathbb{V},\bbSigma )
=E( n^{-1}\widetilde{\mathbb{X}}^{\T}\mathbb{V}^{-1}\bbSigma %
\mathbb{V}^{-1}\widetilde{\mathbb{X}}) $.

The following result gives the asymptotic distribution of $\widehat
{\bolds{%
\beta}}$ for general working covariance matrices.

\begin{theorem}
\label{THM:betahat-normality}Under Assumptions \textup{(A1)}--\textup{(A5)}, as $%
n\rightarrow\infty$,
\[
n^{1/2}( \widehat{\bolds{\beta}}-\bolds{\beta}) \rightarrow
\mathrm{N}( 0,\bolds{\Omega}( \mathbb{V},\bbSigma )
) .
\]
\end{theorem}

\begin{remark} It is easy to show that the covariance
$\bolds{\Omega}%
( \mathbb{V},\bbSigma ) $ in (\ref{EQ:IV}) is minimized by
$\mathbb{V}=\bbSigma $, and in this case equals to $\{
\widetilde{\mathbf{A}}( \mathbb{V}) \} ^{-1}$. To construct
the confidence sets for $\bolds{\beta}$, $\bolds{\Omega}( \mathbb{V}
,\bbSigma ) $ is consistently estimated by
\[
\widehat{\bolds{\Omega}}( \mathbb{V},\widehat{\bbSigma }%
) =n( \widehat{\mathbb{X}}^{\T}\mathbb{V}^{-1}\widehat{\mathbb{X}}%
) ^{-1}( \widehat{\mathbb{X}}^{\T}\mathbb{V}^{-1}\widehat{\bbSigma}\mathbb{V}^{-1}\widehat{\mathbb{X}}) ( \widehat{\mathbb{X%
}}^{\T}\mathbb{V}^{-1}\widehat{\mathbb{X}}) ^{-1},
\]
where $\widehat{\mathbb{X}}=\{ ( \underline{\widehat{\mathbf{X}}}%
_{1}^{\T},\ldots,\underline{\widehat{\mathbf{X}}}_{n}^{\T}) ^{\T%
}\} _{n_{\mathrm{T}}\times d_{1}}$, and
\begin{equation}
\widehat{\underline{\mathbf{X}}}_{i}=\underline{\mathbf{X}}_{i}-\operatorname{Proj} _{G_{n}^{\ast}}\underline{\mathbf{X}}_{i},\qquad i=1,\ldots,n,
\label{EQ:Xhat}
\end{equation}
in which $\operatorname{Proj} _{G_{n}^{\ast}}$ is the projection onto
the empirically centered spline space.
\end{remark}

\begin{remark} The result of Proposition 2 with identity
link in Wang, Carroll and Lin~\cite{WCL05} is a special case of Theorem~\ref{THM:betahat-normality} with $%
\mathbf{V}_{1}=\cdots=\mathbf{V}_{n}=\mathbf{V}$, $m_{1}=\cdots
=m_{n}=M$ and $d_{2}=1$.
\end{remark}

The next theorem shows that the estimated function $\widehat{\eta}%
_{l}$ in (\ref{EQ:mhat}) is $L_{2}$-consistent.

\begin{theorem}
\label{THM:asymptoticphihat}Under Assumptions \textup{(A1)--(A5)}, $\vert
\widehat{\eta}_{l}-\eta
_{l}\vert_{L_{2}}^{2}=\mathrm{O}_{P}\{ J_{n}^{1-2p}+(
J_{n}/n)\}$, for $1\leq l\leq
d_{2}$.
\end{theorem}

\subsection{Sampling properties for the penalized estimators}
\label{SEC:penalizedestimators}

We next show that with a proper choice of $\lambda_{k}$, the penalized
estimator $\widehat{\bolds{\beta}}^{\P}$ has an oracle property. To
avoid confusion, let $\bolds{\beta}_{0}$ be the true value of $\bolds
{\beta}$. Let $r$
be the number of non-zero components of $\bolds{\beta}_{0}$.
Let $\bolds{\beta}_{0}=(\beta_{10},\ldots,\beta_{d_{1}0})^{\T}=(\bolds{\beta}%
_{10}^{\T},\bolds{\beta}_{20}^{\T})^{\T}$, where $\bolds{\beta
}_{10}$ is
assumed to consist of all $r$ non-zero components of $\bolds{\beta}_{0}$,
and $\bolds{\beta}_{20}=\mathbf{0}$ without loss of generality. In a
similar fashion to $\bolds{\beta}$, we can write the collections of all
parametric components, $\mathbb{X}=(\mathbb{X}_{1}^{\T},\mathbb
{X}_{2}^{\T%
})^{\T}$, $\widetilde{\mathbb{X}}=(\widetilde{\mathbb{X}}_{1}^{\T},%
\widetilde{\mathbb{X}}_{2}^{\T})^{\T}$. Denote $a_{n}=\max_{1\leq k\leq
d_{1}}\{|p_{\lambda_{k}}^{\prime}(|\beta
_{k0}|)|,\beta_{k0}\neq0\}$, $w_{n}=\max_{1\leq k\leq
d_{1}}\{|p_{\lambda_{k}}^{\prime\prime}(|\beta_{k0}|)|,\beta
_{k0}\neq
0\}$.

\begin{theorem}
\label{THM:rootn} Under Assumptions \textup{(A1)--(A5)}, and if $a_{n}\to0$ and
$%
w_{n}\to0$ as $n\to\infty$, then there exists a local solution
$\widehat{%
\bolds{\beta}}^{\P}$ in (\ref{DEF:beta-PGEE}) such that its rate of
convergence is $\mathrm{O}_{P}(n^{-1/2}+a_{n})$.
\end{theorem}

Next define a vector $\bolds{\kappa}_{n}=\{p_{\lambda_{1}}^{\prime
}(|\beta_{10}|)\sgn(\beta_{10}),\ldots,p_{\lambda_{r}}^{\prime
}(|\beta_{r0}|)\sgn%
(\beta_{r0})\}^{\T}$ and a diagonal matrix $\bolds{\Sigma}_{\lambda
}=\diag%
\{p_{\lambda_{1}}^{\prime\prime}(|\beta_{10}|),\ldots,p_{\lambda
_{r}}^{\prime\prime}(|\beta_{r0}|)\}$. We further denote $\bolds
{\Sigma}_{1i}=\mathrm{Var}(\underline{\mathbf{Y}}_{i}\vert\underline
{%
\mathbf{X}}_{1i},\underline{\mathbf{Z}}_{i}) $, $\bbSigma _{1}=%
\operatorname{diag}( \bolds{\Sigma} _{11},\ldots,\bolds{\Sigma} _{1n})
$, $\widetilde{\mathbf{A}}_{1}( \mathbb{V}) =E( \widetilde{%
\mathbb{X}}_{1}^{\T}\mathbb{V}^{-1}\widetilde{\mathbb{X}}_{1}) $ and $%
\widetilde{\mathbf{B}}_{1}( \mathbb{V},\bbSigma  _{1} ) =E(
\widetilde{\mathbb{X}}_{1}^{\T}\mathbb{V}^{-1}\bbSigma
_{1}\mathbb{V}%
^{-1}\widetilde{\mathbb{X}}_{1}) $.

The theorem below shows that under regularity conditions, all the covariates
with zero coefficients can be detected simultaneously with probability
tending to 1, and the estimators of all the non-zero coefficients are
asymptotically normally distributed.

\begin{theorem}
\label{THM:oracle} Under Assumptions \textup{(A1)--(A5)}, if
$\lim_{n\rightarrow\infty}\sqrt{n}\lambda_{kn}\rightarrow\infty$ and
\[
\liminf_{n\rightarrow\infty}\liminf_{\beta_{k}\rightarrow
0^{+}}\lambda
_{kn}^{-1}p_{\lambda_{kn}}^{\prime}(|\beta_{k}|)>0,
\]
then the $\sqrt
{n}$-consistent estimator $\widehat{\bolds{\beta}}^{\P}$
in Theorem~\ref{THM:rootn} satisfies $P(\widehat{\bolds{\beta
}}_{2}^{\P}=%
\mathbf{0})\rightarrow1$, as $n\rightarrow\infty$, and
\[
\sqrt{n}\{ \widetilde{\mathbf{A}}_{1}( \mathbb{V}) +\bolds{%
\Sigma} _{\lambda}\} [ \widehat{\bolds{\beta}}_{1}^{{%
\mbox{\rm\tiny P}}}-\bolds{\beta}_{10}+\{ \widetilde{\mathbf{A}}%
_{1}( \mathbb{V}) +\bolds{\Sigma} _{\lambda}\} ^{-1}%
\bolds{\kappa}_{n}] \rightarrow\mathrm{N}(\mathbf{0},\widetilde{%
\mathbf{B}}_{1}( \mathbb{V},\bbSigma _{1} ) ).
\]
\end{theorem}

\section{Implementation}
\label{SEC:Implementation}

In this section, we illustrate how to implement the proposed method in
the semiparametric marginal estimation and variable selection. Let
\[
\bolds{\Sigma}_{\bolds{\lambda}}( \bolds{\beta}) =\operatorname{diag}\biggl\{ \frac{p_{\lambda_{1}}^{\prime}( \vert\beta
_{1}\vert) }{\epsilon+\vert\beta_{1}\vert},\ldots,%
\frac{p_{\lambda_{d_{1}}}^{\prime}( \vert\beta
_{d_{1}}\vert) }{\epsilon+\vert\beta_{d_{1}}\vert
}\biggr\} \label{EQ:Sigma_beta}
\]
for a small number $\epsilon$ ($\epsilon=10^{-6}$ in our simulation
studies). Applying the usual Taylor approximation, $Q_{\mathcal{P}}(
\bolds{\beta}) $ can be locally approximated by
\[
Q( \bolds{\beta}) +\dot{Q}( \bolds{\beta}_{0})
^{T}( \bolds{\beta-\beta}_{0}) +\tfrac{1}{2}( \bolds{%
\beta-\beta}_{0}) ^{T}\ddot{Q}( \bolds{\beta}_{0})
( \bolds{\beta-\beta}_{0}) +\tfrac{1}{2}n_{\T}\bolds{\beta}%
^{T}\Sigma_{\bolds{\lambda}}( \bolds{\beta}_{0}) \bolds{%
\beta}.
\]
By the local quadratic approximations for penalty functions (Fan and Li
\cite{FL01}, Section 3.3), the solution can be found iteratively,
\[
\bolds{\beta}^{( k+1) }=\Biggl[ \sum_{i=1}^{n}\widehat{%
\underline{\mathbf{X}}}_{i}^{\T}\bigl\{ \mathbf{V}_{i}^{( k)
}\bigr\} ^{-1}\widehat{\underline{\mathbf{X}}}_{i}+n_{\T}\Sigma_{\bolds{%
\lambda}}\{ \bolds{\beta}^{( k) }\} \Biggr]
^{-1}\sum_{i=1}^{n}\widehat{\underline{\mathbf{X}}}_{i}^{\T}\bigl\{ \mathbf{V
}_{i}^{( k) }\bigr\} ^{-1}\{ \underline{\mathbf{Y}}_{i}-%
\widehat{\Pi}_{n}\underline{\mathbf{Y}}_{i}\} ,
\]
where $\widehat{\Pi}_{n}\underline{\mathbf{Y}}_{i}$ is the projection
of $\underline{\mathbf{Y}}_{i}$ onto the spline space $G_{n}^{0}$, and $
\underline{\widehat{\mathbf{X}}}_{i}$ is given in (\ref{EQ:Xhat}).

Following Fan and Li~\cite{FL01}, we derive a sandwich formula for the
standard errors
of the estimated covariates $\widehat{ \bolds{\beta}}^{{\P}}$%
\begin{eqnarray} \label{EQ:std}
\widehat{\operatorname{Cov}}( \widehat{\bolds{\beta}}^{{\P}})
=\{ \ddot{Q}( \widehat{\bolds{\beta}}_{\bolds{\lambda}}^{{\P
}}) +n_{\mathrm{T}}\Sigma_{\bolds{\lambda}}( \widehat{\bolds{%
\beta}}_{\bolds{\lambda}}^{{\P}}) \} ^{-1}\widehat{\mathrm{%
Cov}}\{ \dot{Q}( \widehat{\bolds{\beta}}_{\bolds{\lambda}}^{{%
\P}}) \}
\times\{ \ddot{Q}( \widehat{\bolds{\beta}}_{\bolds{\lambda}%
}^{{\P}}) +n_{\mathrm{T}}\Sigma_{\bolds{\lambda}}( \widehat{%
\bolds{\beta}}_{\bolds{\lambda}}^{{\P}}) \} ^{-1},
\end{eqnarray}
where $\ddot{Q}( \bolds{\beta} ) =\sum_{i=1}^{n}\widehat{\underline{%
\mathbf{X}}}_{i}^{\T}\mathbf{V}_{i}^{-1}\widehat{\underline{\mathbf{X}}}_{i}$
and $\widehat{\operatorname{Cov}}\{ \dot{Q}( \bolds{\beta} ) \}
=\sum_{i=1}^{n}\widehat{\underline{\mathbf{X}}}_{i}^{\T}\mathbf
{V}_{i}^{-1}%
\widehat{\bolds{\Sigma}}_{i}\mathbf{V}_{i}^{-1}\widehat{\underline{\mathbf{X}}}_{i}$. Applying conventional techniques that arise in the likelihood
setting, we can show that the above sandwich formula is a consistent estimator
and has good accuracy in our simulation study for moderate sample sizes.

We use BIC to select the tuning parameters $\bolds{\lambda}=(
\lambda_{1},\ldots,\lambda_{d_{1}}) $. Let
\[
e( \bolds{\lambda}) =\mathrm{tr}\{ [ \ddot{Q}(
\widehat{\bolds{\beta}}_{\bolds{\lambda}}^{\P}) +n_{\mathrm{T}%
}\Sigma_{\bolds{\lambda}}( \widehat{\bolds{\beta}}_{\bolds{%
\lambda}}^{\P}) ] ^{-1}\ddot{Q}( \widehat{\bolds{\beta}}%
_{\bolds{\lambda}}^{{\P}}) \}
\]
be the effective number of parameters in the last step of the Newton--Raphson
iteration. Then
\[
\mathrm{BIC}( \bolds{\lambda}) =\log\Biggl\{ \frac{1}{n_{\mathrm{T}
}}\sum_{i=1}^{n}(\mathbf{y}_{i}-\widehat{\boldmu}_{i})^{\T}\mathbf{R}_i^{-1}
(\mathbf{y}_{i}-\widehat{\boldmu}_{i})\Biggr\} +\frac{\log
( n_{\mathrm{T}})}{n_{\mathrm{T}}} e( \bolds{\lambda}).
\]
The minimization problem over a $d$-dimensional space is difficult. However,
Li and Liang~\cite{LL08} conjectured that the magnitude of $\lambda_{k}$ should be
proportional to the standard error of $\beta_{k}$. So we suggest
taking $%
\lambda_{k}=\lambda*\mathrm{SE}( \widehat{\beta}_{k}) $, in
practice, where $\mathrm{SE}( \widehat{\beta}_{k}) $ is the
standard error of $\widehat{\beta}_{k}$, the unpenalized estimator defined
above. Thus, the minimization problem can be reduced to a one-dimensional
problem, and the tuning parameter can be estimated by a grid search.

\section{Simulation}
\label{SEC:Simulation}

In this section, we discuss finite sample properties of the proposed
estimators via simulation studies. We simulated $100$ data sets of size
$%
n=100$, $200$ and $400$ from the model
\begin{equation} \label{SIM:model}
Y_{ij}=\bolds{\beta^T} \mathbf{X}_{ij}+\eta_1(Z_{ij1})+\eta_2(Z_{ij2})+
\varepsilon_{ij},\qquad  i=1,\ldots,n, j=1,\ldots,3
\end{equation}
where the coefficients $\bolds{\beta}=(3, 1.5, 0, 0, 2, 0, 0, 0)^{\T}$,
function $\eta_1(z) = \sin{2\uppi(z-0.5)}$ and function $\eta_2(z) =
z-0.5+\sin\{2\uppi(z-0.5)\}$.

The $2$-vector $\underline{\mathbf{Z}}_i$ was generated from a
bivariate normal distribution with mean $0$, a common marginal variance
$0.25$ with correlation $ 0.9$, but truncated to the unit square $%
[0, 1]^2$. The covariates $X_{ijk}$, $k=1,\ldots,6$, were generated
independently from N$(0,0.25)$. Covariate $X_{ij7} =
3(1-2Z_{ij1})(1-2Z_{ij2}) +
u_{ij}$, where $u_{ij}\sim N(0, 0.25)$ and is independent of $\mathbf
{Z}_{ij}$. Covariate $X_{ij8}$ was generated as $-0.5$ and $0.5$ with
equal probability. We generated $\bolds{\varepsilon}_i=(\varepsilon_{i1},
\varepsilon_{i2}, \varepsilon_{i3})$ from N$(0,\Sigma_E)$, where $\Sigma
_E=(1-\alpha)\mathbf{I} +\alpha\mathbf{11}^{\T}$ with $\mathbf{1}$
being a vector with all ``$1$'' and $\alpha=0.9$, that is, $\Sigma_E$
is exchangeable.

Cubic B-splines were used to approximate the nonparametric functions as
described in Section~\ref{SEC:Methodology}. We tried different numbers of knots (ranging
from $2$ to $10$) and found that the choice of number of knots didn't
make a significant difference in this simulation study. Our reported
results in Tables~\ref{TAB:selection} and~\ref{TAB:betaerror} were
based on
using $4$ equally spaced knots.

%
\begin{table}
\tabcolsep=0pt
\caption{Model selection and estimation: the average number of correct
(C) and incorrect (I) $0$s, MRME (\%) and RMSE}
\label{TAB:selection}
\begin{tabular*}{\textwidth}{@{\extracolsep{\fill}}llllllllllllll@{}}
\hline
& & \multicolumn{4}{l}{EX} & \multicolumn{4}{l}{AR(1)} & \multicolumn{4}{l@{}}{WI}
\\[-6pt]
& & \multicolumn{4}{c}{\hrulefill} & \multicolumn{4}{c}{\hrulefill} & \multicolumn{4}{c@{}}{\hrulefill}
\\
$n$& Penalty& C & I & MRME & RMSE & C & I & MRME & RMSE & C & I & MRME & RMSE\\
\hline
$100$ & SCAD & 4.67 & 0 & 80.63 &0.1592 & 4.64 & 0 & 84.65 &0.1727 & 4.64 &
0 & 82.38 &0.5883\\
& HARD & 4.80 & 0 & 85.90 &0.1691 & 4.70 & 0 & 86.56 & 0.1916& 4.85 &
0 & 85.81 &0.4410\\
& ORACLE & 5.00 & 0 & 77.23 &0.1586& 5.00 & 0 & 73.40 &0.1723 & 5.00 &
0 & 70.71 &0.4126\\[5pt]
$200$& SCAD & 4.72 & 0 & 76.30 &0.1053 & 4.72 & 0 & 81.63 &0.1127 & 4.70 &
0 & 79.99 &0.3921\\
& HARD & 4.79 & 0 & 82.81 &0.1116 & 4.71 & 0 & 82.18 &0.1252 & 4.98 &
0 & 86.15 &0.2816\\
& ORACLE & 5.00 & 0 & 66.96 &0.1038 & 5.00 & 0 & 66.18 &0.1110 & 5.00
& 0 & 70.86 &0.2787\\[5pt]
$400$ & SCAD & 4.92 & 0 & 84.91 &0.0733 & 4.86 & 0 & 84.50 &0.0864 & 4.88 &
0 & 85.78 &0.2689\\
& HARD & 4.93 & 0 & 91.23 &0.0758 & 4.87 & 0 & 85.65 &0.0924 & 4.90 &
0 & 91.08 &0.2021\\
& ORACLE & 5.00 & 0 & 68.33 &0.0731 & 5.00 & 0 & 66.91 &0.0860 & 5.00
& 0 & 71.52 &0.1857\\
\hline
\end{tabular*}
\end{table}

To the simulated data sets, we applied the proposed method for
estimation and variable selection. To study how the structure of the
working correlation could affect our estimation and variable selection
results, we considered the following three correlation structures: the
correct exchangeable working correlation structure (EX), working
independence (WI) and AR (1) structures. Table~\ref{TAB:selection}
summarizes the estimation and variable selection results with two types
of penalty functions: SCAD and HARD. The average number of zero
coefficients is reported in Table~\ref{TAB:selection}, in which the
column labeled ``C'' presents the average restricted only to the true
zero coefficients, and the column labeled ``I'' shows the average of
numbers erroneously set to zero. The rows with ``SCAD'' and ``HARD''
stand, respectively, for the penalized least squares with the SCAD and
HARD penalties. The oracle estimates always identify the 5 zero
coefficients and 3 non-zero coefficients correctly. The medians of
relative model errors (MRME) as suggested in Fan and Li~\cite{FL01} and
the root mean squared errors (RMSE) of the estimated coefficients over
100 simulated data sets are also reported in Table~\ref{TAB:selection}.

From Table~\ref{TAB:selection}, one sees that the choice of correlation
structure has little impact on the results of variable selection: the
number of correctly identified zero coefficients are all close to 5
regardless the correlation structure; and none of the nonzero
coefficients were erroneously set to 0 in any scenario. Table \ref
{TAB:selection} also shows that the estimators with correct working
correlation have the smallest RMSEs, thus are more efficient than those
estimators with misspecified working correlation structures. The
efficiency of the estimators based on the AR(1) is close to those based
on EX, but there seems to be some significant loss of efficiency for
the estimators based on the WI structure which ignores the within
subject/cluster correlation. In terms of choosing penalty functions, we
find that both HARD and SCAD perform very well and the corresponding
MRME and RMSE are comparable to those of the ORACLE.

%
\begin{table}
\tabcolsep=0pt
\caption{Simulation results on standard error estimation for the non-zero
coefficients ($\protect\beta_1, \protect\beta_2, \protect\beta_5$)}
\label{TAB:betaerror}
\begin{tabular*}{\textwidth}{@{\extracolsep{\fill}}lllllllllll@{}}
\hline
 &  & \multicolumn{3}{l}{$\hat{\beta_1}$} & \multicolumn{3}{l}{$\hat{\beta_2}$} & \multicolumn{3}{l@{}}{$\hat{\beta_5}$}
 \\[-6pt]
  &  & \multicolumn{3}{c}{\hrulefill} & \multicolumn{3}{c}{\hrulefill} & \multicolumn{3}{c@{}}{\hrulefill}
 \\
$n$& Penalty& SD & SD$_\mathrm{m}$ & SD$_{\mathrm{mad}}$ & SD & SD$_\mathrm{m}$ & SD$_{\mathrm{mad}}$ & SD & SD$_\mathrm{m}$
& SD$%
_{\mathrm{mad}}$ \\
\hline
EX & & & & & & & & & \\
$100$
& SCAD & 0.0889 & 0.0906 & 0.0141 & 0.1082 & 0.0911 & 0.0116 & 0.1034 &
0.0894 & 0.0111\\
& HARD & 0.0879 & 0.0907 & 0.0118 & 0.1102 & 0.0911 & 0.0113 & 0.0982 &
0.0897 & 0.0108\\
& ORACLE & 0.0866 & 0.0988 & 0.0062 & 0.1066 & 0.0899 & 0.0112 & 0.1012
& 0.0903 & 0.0088
\\ [3pt]
$200$& SCAD & 0.0655 & 0.0638 & 0.0035 & 0.0616 & 0.0629 & 0.0036 & 0.0594 &
0.0633 & 0.0036\\
& HARD & 0.0655 & 0.0637 & 0.0033 & 0.0627 & 0.0630 & 0.0033 & 0.0600 &
0.0632 & 0.0035\\
& ORACLE & 0.0648 & 0.0699 & 0.0078 & 0.0614 & 0.0637 & 0.0034 & 0.0594
& 0.0629 & 0.0036
\\ [3pt]
$400$
& SCAD & 0.0414 & 0.0445 & 0.0043 & 0.0379 & 0.0445 & 0.0041 & 0.0415 &
0.0449 & 0.0042\\
& HARD & 0.0414 & 0.0446 & 0.0043 & 0.0373 & 0.0445 & 0.0041 & 0.0418 &
0.0449 & 0.0042\\
& ORACLE & 0.0412 & 0.0485 & 0.0086 & 0.0368 & 0.0443 & 0.0049 & 0.0404
& 0.0445 & 0.0046
\\ [5pt]
AR(1) & & & & & & & & & \\
$100$
& SCAD & 0.0983 & 0.0923 & 0.0141 & 0.1035 & 0.0940 & 0.0129 & 0.1153 &
0.0924 & 0.0132\\
& HARD & 0.0996 & 0.0920 & 0.0160 & 0.1073 & 0.0939 & 0.0130 & 0.1117 &
0.0924 & 0.0128\\
& ORACLE & 0.0976 & 0.0972 & 0.0097 & 0.0971 & 0.0915 & 0.0124 & 0.1173
& 0.0930 & 0.0122
\\ [3pt]
$200$
& SCAD & 0.0635 & 0.0646 & 0.0041 & 0.0539 & 0.0634 & 0.0047 & 0.0647 &
0.0639 & 0.0045\\
& HARD & 0.0632 & 0.0645 & 0.0045 & 0.0544 & 0.0634 & 0.0044 & 0.0626 &
0.0639 & 0.0045\\
& ORACLE & 0.0624 & 0.0689 & 0.0073 & 0.0535 & 0.0646 & 0.0049 & 0.0657
& 0.0635 & 0.0048
\\[3pt]
$400$
& SCAD & 0.0452 & 0.0448 & 0.0056 & 0.0390 & 0.0451 & 0.0055 & 0.0535 &
0.0452 & 0.0052\\
& HARD & 0.0451 & 0.0449 & 0.0057 & 0.0390 & 0.0451 & 0.0055 & 0.0537 &
0.0452 & 0.0053\\
& ORACLE & 0.0454 & 0.0477 & 0.0061 & 0.0392 & 0.0448 & 0.0051 & 0.0539
& 0.0450 & 0.0056
\\[5pt]
WI & & & & & & & & & \\
$100$
& SCAD & 0.2177 & 0.2364 & 0.0164 & 0.2192 & 0.2381 & 0.0187 & 0.2341 &
0.2375 & 0.0189\\
& HARD & 0.2235 & 0.2364 & 0.0151 & 0.2185 & 0.2396 & 0.0169 & 0.2393 &
0.2389 & 0.0193\\
& ORACLE & 0.2239 & 0.0579 & 0.1623 & 0.2118 & 0.2341 & 0.0165 & 0.2218
& 0.2374 & 0.0181
\\[3pt]
$200$
& SCAD & 0.1864 & 0.1674 & 0.0204 & 0.1697 & 0.1677 & 0 .0188 & 0.1328
& 0.1671 & 0.0186\\
& HARD & 0.1876 & 0.1675 & 0.0201 & 0.1676 & 0.1680 & 0.0180 & 0.1385 &
0.1676 & 0.0181\\
& ORACLE & 0.1836 & 0.0415 & 0.1242 & 0.1656 & 0.1669 & 0.0162 & 0.1356
& 0.1671 & 0.0165
\\[3pt]
$400$& SCAD & 0.0957 & 0.1162 & 0.0131 & 0.1055 & 0.1163 & 0.0117 & 0.1252 &
0.1164 & 0.0131\\
& HARD & 0.0956 & 0.1162 & 0.0129 & 0.1042 & 0.1165 & 0.0121 & 0.1222 &
0.1163 & 0.0128\\
& ORACLE & 0.0956 & 0.0289 & 0.0778 & 0.1069 & 0.1160 & 0.0110 & 0.1227
& 0.1162 & 0.0104
\\
\hline
\end{tabular*}
\end{table}

We also tested the accuracy of our standard error formula based on (\ref
{EQ:std}). The median absolute deviation (MAD) divided by 0.6745
(denoted by
$\mathrm{SD}$ in Table~\ref{TAB:betaerror}) of 100 estimated
coefficients from
the 100 simulations can be regarded as the true standard error. The median
of the 100 estimated SDs (denoted by $\mathrm{SD}_\mathrm{m}$) and the
MAD error of the 100 estimated standard errors divided by 0.6745
(denoted by $\mathrm{SD}_{\mathrm{mad}}$) gauge the overall performance
of the standard error. Table~\ref{TAB:betaerror} presents the standard
errors for non-zero coefficients when the sample size $n =100$, $200$
and $400$. It suggests
that the sandwich formula performs satisfactorily for SCAD and HARD penalties.
The standard errors based on the SCAD and HARD penalty functions are
closer to those of the ORACLE as $n$ increases. Similarly to the RMSE
results shown in Table~\ref{TAB:selection}, Table~\ref{TAB:betaerror}
also shows that the estimation procedures with a correct EX working
correlation are more efficient than their counterparts with WI working
correlation. Estimation based on a misspecified AR(1) correlation
structure will lead to some efficiency loss, but it is quite close to
using the true EX structure.\vadjust{\goodbreak}

\section{Application}
\label{SEC:Application}

To illustrate our method, we considered the longitudinal CD4 cell count
data among HIV seroconverters. This dataset contains $2376$
observations of CD4 cell counts on $369$ men infected with the HIV virus;
see Zeger and Diggle~\cite{ZD94} for a detailed description of this dataset. Both Wang, Carroll and Lin~\cite{WCL05}
and Huang, Zhang and Zhou~\cite{HZZ06} analyzed the same dataset using a PLM.
Their analysis aimed to estimate the average time course of CD4 counts and
the effects of other covariates. In our analysis, we fit the data using an
APLM, with the square root transformed CD4 counts
as the response, and covariates including AGE, SMOKE (smoking status
measured by packs of cigarettes), DRUG (yes, 1; no, 0), SEXP (number of sex
partners), DEPRESSION (measured by the CESD scale) and YEAR (the effect of
time since seroconversion). To take advantage of flexibility of partially
linear additive models, we let both DEPRESSION and YEAR be modeled
nonparametrically, the remaining parametrically. It is of interest to
examine whether there are any interaction effects between the parametric
covariates, so we included all these interactions in the parametric part.

%
\begin{sidewaystable}
\tablewidth=\textwidth
\caption{Estimated coefficients for CD4 dataset}
\label{TAB:APPest}
\begin{tabular*}{\textwidth}{@{\extracolsep{\fill}}ld{2.10}d{2.10}d{2.10}d{2.10}d{2.10}d{2.10}@{}}
\hline
 & \multicolumn{3}{l}{Full} & \multicolumn{3}{l@{}}{Penalized}  \\[-6pt]
  & \multicolumn{3}{c}{\hrulefill} & \multicolumn{3}{c@{}}{\hrulefill}  \\
\multicolumn{1}{@{}l}{Variable}&
\multicolumn{1}{l}{WI} &
\multicolumn{1}{l}{AR(1)}&
\multicolumn{1}{l}{RSM} &
\multicolumn{1}{l}{WI} &
\multicolumn{1}{l}{AR(1)} &
\multicolumn{1}{l}{RSM} \\
& \multicolumn{1}{l}{$\hat{\beta}$ ($\mathrm{SE}(\hat{\beta})$)}&
\multicolumn{1}{l}{$\hat{\beta}$ ($\mathrm{SE}(\hat{\beta})$)} &
\multicolumn{1}{l}{$\hat{\beta}$ ($\mathrm{SE}(\hat{\beta})$)} &
\multicolumn{1}{l}{$\hat{\beta}$ ($\mathrm{SE}(\hat{\beta})$)} &
\multicolumn{1}{l}{$\hat{\beta}$ ($\mathrm{SE}(\hat{\beta})$)} &
\multicolumn{1}{l@{}}{$\hat{\beta}$ ($\mathrm{SE}(\hat{\beta})$)} \\
\hline
INTERCEPT &  24.365\ (0.417)  &  24.540\ (0.480)  &  24.819\ (0.494)  &
 24.487\ (0.391)   & 24.454\ (0.464) & 24.793\ (0.461) \\
AGE &   -0.013\ (0.035)  &  -0.023\ (0.045)  &  -0.049\ (0.049)  &
 \multicolumn{1}{l}{\phantom{0}0 (0)}  & \multicolumn{1}{l}{\phantom{0}0 (0)}  & \multicolumn{1}{l}{\phantom{0}0 (0)}   \\
SMOKE &   1.070\ (0.234)  &  0.825\ (0.259)  &  0.654\ (0.264)  &
0.733\ (0.148)  &  0.824\ (0.247) &  0.424\ (0.176) \\
DRUG &   2.671\ (0.491)  &  1.958\ (0.517)   &  1.468\ (0.507)  &   2.486\ (0.454)  &  2.025\ (0.511) & 1.340\ (0.462)   \\
SEXP &   0.165\ (0.082)  &  0.153\ (0.084)  &  0.109\ (0.082)  &
0.170\ (0.076)   &  0.174\ (0.080) &  0.144\ (0.078)  \\
AGE*SMOKE &   -0.014\ (0.011)  &   0.000\ (0.015)  &  0.002\ (0.017)
&  \multicolumn{1}{l}{\phantom{0}0 (0)}  & \multicolumn{1}{l}{\phantom{0}0 (0)}  & \multicolumn{1}{l}{\phantom{0}0 (0)}  \\
AGE*DRUG &  0.043\ (0.036)  &   0.008\ (0.042)  &  0.010\ (0.044)  &
 \multicolumn{1}{l}{\phantom{0}0 (0)}  &  \multicolumn{1}{l}{\phantom{0}0 (0)}  & \multicolumn{1}{l}{\phantom{0}0 (0)}  \\
AGE*SEXP &   0.001\ (0.005)  &  0.005\ (0.005)  &  0.008\ (0.005)  &
 \multicolumn{1}{l}{\phantom{0}0 (0)}  & \multicolumn{1}{l}{\phantom{0}0 (0)}  & \multicolumn{1}{l}{\phantom{0}0 (0)}  \\
SMOKE*DRUG &   -0.402\ (0.233)  &  -0.331\ (0.246)  &   -0.288\ (0.248)  &  \multicolumn{1}{l}{\phantom{0}0 (0)}   &  -0.337\ (0.241) & \multicolumn{1}{l}{\phantom{0}0 (0)}   \\
SMOKE*SEXP &   0.058\ (0.024)  &   0.043\ (0.026)   &  0.047\ (0.025)  &
 0.051\ (0.023)  &  0.045\ (0.025) &  0.046\ (0.025)  \\
DRUG*SEXP &   -0.364\ (0.087)  &   -0.251\ (0.086)  &  -0.17\ (0.083)  &
 -0.355\ (0.084)  &  -0.265\ (0.084) &  -0.186\ (0.081) \\
 \hline
\end{tabular*}
\end{sidewaystable}

For the working variance, we considered the WI, the AR(1) and the
``random intercept plus serial correlation and measurement error''
covariance (RSM) in Zeger and Diggle~\cite{ZD94}. One can obtain the RSM
structure by fitting a full model to the data and inspecting the
variogram of the residuals. Wang, Carroll and Lin~\cite{WCL05}
and Huang, Zhang and Zhou~\cite{HZZ06} also analyzed this data set using the RSM structure.
More precisely, the working covariance matrices are specified by $\tau
^2 \mathbf{I} +\nu^2 \mathbf{J} +\omega^2 \mathbf{H}$, where $\mathbf
{I}$ is an identity matrix, $\mathbf{J}$ is a matrix of 1s and $\mathbf
{H}(j,j^\prime)=\exp(-\alpha|\mathrm{YEAR}_{ij} -\mathrm
{YEAR}_{ij^\prime} |)$. We used the covariance parameters $(\tau^2,\nu
^2,\omega^2,\alpha^2)=(11.32,3.26,22.15,0.23)$ calculated by Wang
\textit{et al.}~\cite{WCL05}. Table~\ref{TAB:APPest} gives the estimates of
the regression coefficients using WI, AR(1) and RSM covariance structures.
The standard errors (SE) were all calculated using the sandwich method.
We used cubic splines of $4$ knots selected by the five-fold
delete-subject-out cross-validation from the range of 0--20. We refer
the reader to Huang, Wu, and Zhou~\cite{HWZ04} for the detail of the
delete-subjects-out $K$-fold cross-validation. The left panel of Table
\ref{TAB:APPest} reports the estimation using full model, and the
selection results are shown in the right panel.

We further applied the proposed approach to select significant
variables. We used the SCAD penalty, the tuning parameter $\lambda
=0.4549, 0.2829, 0.3143$ for WI, AR(1) and RSM covariance structure,
respectively. The results are also shown in Table~\ref{TAB:APPest}.
Under both WI and RSM structures, SMOKE, DRUGS, SEXP, SOMKE$\ast$SEXP and
DRUGS$\ast$SEXP are identifies as significant covariates. One notes some
slight selection difference when AR(1) structure is used, which
suggests that SMOKE$\ast$DRUGS may also be significant. Although the
selection procedure is not sensitive to the choice of covariance
structure as shown in our simulation study, different covariance
structures may still lead to slight different results. Therefore, it is
important for one to choose a
covariance structure close to the true one. We also find some
significant interactions among some covariates which may be ignored by
Wang, Carroll and Lin~\cite{WCL05} and Huang, Zhang and Zhou~\cite{HZZ06}.

\begin{figure}

\includegraphics{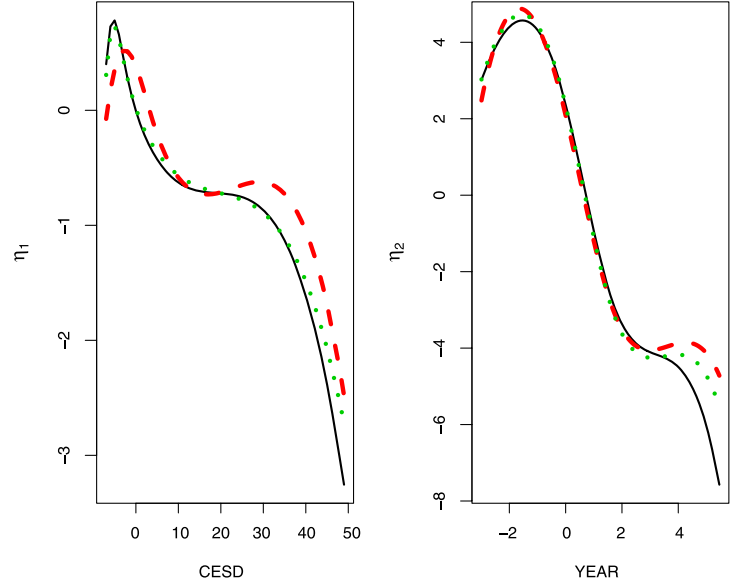}

\caption{The estimates of the nonpararmetric components: $\widehat
{\eta}_1$ and $\widehat{\eta}_2$. The solid, dotted and dashed curves
correspond to the estimates under WI, AR(1) and RSM structures.}
\label{FIG:cd4}
\end{figure}

The nonparametric curve estimates using the WI (solid line), AR(1)
(dotted line) and RSM (dashed line) estimators are plotted in Figure
\ref{FIG:cd4} for
``DEPRESSION'' and ``YEAR.'' One can see that it is more reasonable to put
``DEPRESSION'' as a nonparametric component.

\section{Discussion}
\label{SEC:Discussion}

We have developed a general methodology for simultaneously selecting
variables and estimating the unknown components in APLMs for
longitudinal and clustered data. We propose a one-step least squares
approach to obtain the estimation of both the parametric and
nonparametric components based on polynomial spline smoothing. This
approach is flexible, computationally simple and very easy to implement
in practice. We demonstrate that the asymptotic normality of the
estimated coefficients for the linear part is retained. The proposed
penalized regression method also achieves an ``oracle'' property in the
sense that it performs as well as if the subset of significant
parametric components were known in advance.

In this paper, our primary interest is the linear components, and we
treat the nonparametric functions as nuisance components; thus we limit
our discussions to estimation and variable selection for the linear
part. Nonetheless, this may be extended to the nonparametric components
using techniques similar to those in Xue~\cite{X09}. An anonymous referee
pointed out the feasibility of obtaining the asymptotic ``oracle''
property of the nonparametric components in Ma and Yang~\cite{MY11}. We believe
that this property can be similarly obtained via a two-step spline
backfitted kernel smoothing procedure (Ma and Yang~\cite{MY11}). However, the
technical details deserve careful consideration, and this is an
interesting topic of future research.\

The simulation result indicates that the variable selection is
consistent even if the correlation structure is misspecified. However,
misspecification may lead to some efficiency loss. So, it would be
desirable if one could choose an appropriate correlation structure
based on available data in practice. The simulation results clearly
show that there is marked improvement of efficiency when one uses the
correct correlation structure though the variable selection seems to be
consistent with misspecified structure. To select the correlation
matrix, one might consider some resampling-based methods, such as the
bootstrap and cross-validation methods in Pan and Connett~\cite{PC02} and other
techniques in Diggle \textit{et al.}~\cite{DHLZ02}. There is, however, a clear
need to formalize the procedures with solid theoretical justification.
Instead of modeling the correlation through the ``working'' correlation
matrix, one could also nonparametrically model the variance--covariance
as some unknown smooth function
(Chiou and M\"{u}ller~\cite{CM05}). This is an excellent research
problem for future study.

\begin{appendix}\label{appendix}

\section*{Appendix}

\renewcommand{\thelemma}{A.\arabic{lemma}}
\renewcommand{\theequation}{A.\arabic{equation}} %
\setcounter{equation}{0}
\setcounter{lemma}{0}

For any vector $\mathbf{x}=(
x_{1},\ldots,x_{d}) ^{\T}$, we denote $\Vert\cdot\Vert$ the
usual Euclidean norm, that is, $\Vert\mathbf{x}\Vert=%
\sqrt{\sum_{k=1}^{d}x_{k}^{2}}$, and $\Vert\cdot\Vert_{\infty}$ the sup
norm, that is, $\Vert\mathbf{x}\Vert_{\infty}=\sup_{1\leq k\leq d}|x_{k}|$.
For any functions $\phi,\varphi$, let $\phi( \underline{\mathbf{X}}%
_{i},\underline{\mathbf{Z}}_{i}) $ and $\varphi( \underline{%
\mathbf{X}}_{i},\underline{\mathbf{Z}}_{i}) $ be $m_{i}$-vectors; then
define the empirical inner product and the empirical norm
as $\langle\phi,\varphi\rangle_{n}\equiv\langle\phi,\varphi
\rangle_{n,\mathbf{V}}=n^{-1}\sum_{i=1}^{n}\phi( \underline{%
\mathbf{X}}_{i},\underline{\mathbf{Z}}_{i}) ^{\T}\mathbf{V}%
_{i}^{-1}\varphi( \underline{\mathbf{X}}_{i},\underline{\mathbf{Z}}%
_{i}) $, $\Vert\phi\Vert
_{n}^{2}=\langle\phi,\phi\rangle_{n} $, for the
working covariance $\mathbf{V}_{i}$. Further denote $E_{n}( \phi
) =\break n^{-1}\sum_{i=1}^{n}\mathbf{1}_{m_{i}}^{\T}\mathbf{V}_{i}^{-1}\phi
( \underline{\mathbf{X}}_{i},\underline{\mathbf{Z}}_{i}) $. If
functions $\phi,\varphi$ are $L^{2}$-integrable, we define the theoretical
inner product and its corresponding theoretical $L^{2}$ norm as $%
\langle\phi,\varphi\rangle=E( \langle\phi
,\varphi\rangle_{n}) $, $\Vert\phi\Vert
^{2}=E( \Vert\phi\Vert_{n}^{2}) $. Let $\widehat{%
\Pi}_{n}$ and $\Pi_{n}$ denote, respectively, the projection onto $G_{n}^{0}$
relative to the empirical and theoretical inner products. For
convenience, let $h=h_{n}\sim J_{n}^{-1}$ and $\mathbf{I}_{d}$ be the
$d\times d$ identity matrix.

\subsection{\texorpdfstring{Proof of Theorem \protect\ref{THM:betahat-normality}}
{Proof of Theorem 1}}\label{SUBSEC:proof2}

\begin{lemma}
\label{LEM:AnBnCn} Define
\begin{eqnarray*}
A_{n}&=&\sup_{g_{1},g_{2}\in G_{n}^{0}}\vert
\langle g_{1},g_{2}\rangle_{n}-\langle
g_{1},g_{2}\rangle\vert\Vert g_{1}\Vert
^{-1}\Vert g_{2}\Vert^{-1}, \\
B_{n}&=&\max_{1\leq k\leq
d_{1}}\sup_{g\in G_{n}^{0}}\bigl\vert\Vert x_{k}-g\Vert
_{n}^{2}/\Vert x_{k}-g\Vert^{2}-1\bigr\vert,
\end{eqnarray*}
then $A_{n}=\mathrm{O}_{P}\{ \sqrt{\log( n) /( nh^{2}) } \}$ and $B_{n}=\mathrm{O}_{P}\{
\sqrt{\log( n) /( nh^{2}) }%
\} $.
\end{lemma}

Lemma~\ref{LEM:AnBnCn} can be proved similarly to Lemmas A2 and A3 in
Huang, Zhang and Zhou~\cite{HZZ06} and are thus omitted.

To obtain the closed-form expression of $\widehat{\bolds{\beta}}$, we need
the following block form of the inverse of $\sum_{i=1}^{n}\underline{\mathbf{D}}_{i}^{\T}\mathbf{V}_{i}^{-1}\underline{\mathbf{D}}_{i} $:
\begin{equation}
\pmatrix{
\displaystyle\sum_{i=1}^{n}\underline{\mathbf{X}}_{i}^{\T}\mathbf{V}_{i}^{-1}%
\underline{\mathbf{X}}_{i} & \displaystyle\sum_{i=1}^{n}\underline{\mathbf{X}}%
_{i}^{\T}\mathbf{V}_{i}^{-1}\underline{\mathbf{B}}_{i} \vspace*{2pt}\cr
\displaystyle\sum_{i=1}^{n}\underline{\mathbf{B}}_{i}^{\T}\mathbf{V}_{i}^{-1}%
\underline{\mathbf{X}}_{i} & \displaystyle\sum_{i=1}^{n}\underline{\mathbf{B}}%
_{i}^{\T}\mathbf{V}_{i}^{-1}\underline{\mathbf{B}}_{i}}^{-1}=\pmatrix{
\mathbf{H}_{\mathbf{XX}} & \mathbf{H}_{\mathbf{XB}} \vspace*{2pt}\cr
\mathbf{H}_{\mathbf{BX}} & \mathbf{H}_{\mathbf{BB}}}^{-1}=\pmatrix{
\mathbf{H}^{11} & \mathbf{H}^{12} \vspace*{2pt}\cr
\mathbf{H}^{21} & \mathbf{H}^{22}} , \label{EQ:DVDinverse}
\end{equation}
where $\mathbf{H}^{11}=( \mathbf{H}_{\mathbf{XX}}-\mathbf{H}_{\mathbf{XB
}}\mathbf{H}_{\mathbf{BB}}^{-1}\mathbf{H}_{\mathbf{BX}}) ^{-1}$, $%
\mathbf{H}^{22}=( \mathbf{H}_{\mathbf{BB}}-\mathbf{H}_{\mathbf{BX}}%
\mathbf{H}_{\mathbf{XX}}^{-1}\mathbf{H}_{\mathbf{XB}}) ^{-1}$, $%
\mathbf{H}^{12}=-\mathbf{H}^{11}\mathbf{H}_{\mathbf{XB}}\mathbf
{H}_{\mathbf{BB}}^{-1}$ and $\mathbf{H}^{21}=-\mathbf{H}^{22}\mathbf{H}_{\mathbf{BX}}
\mathbf{H}_{\mathbf{XX}}^{-1}$. Consequently,
\begin{equation}
\widehat{\bolds{\beta}}=\mathbf{H}^{11}\Biggl\{ \sum_{i=1}^{n}\underline{\mathbf{X}}_{i}^{\T}\mathbf{V}_{i}^{-1}\underline{\mathbf{Y}}_{i}-\mathbf{H}_{\mathbf{XB}}\mathbf{H}_{\mathbf{BB}}^{-1}\sum_{i=1}^{n}\underline{\mathbf{B}}_{i}^{\T}\mathbf{V}_{i}^{-1}\underline{\mathbf{Y}}_{i}\Biggr\}
. \label{EQ:betahat}
\end{equation}

\begin{lemma}
\label{LEM:SandHposdef} Under Assumptions \textup{(A1)--(A5)}, for $\mathbf
{H}_{\mathbf{BB}}$ in \textup{(\ref{EQ:DVDinverse})}, one has \textup{(i)} there exist
constants $0<c_H<C_{H}$, $%
C_{H}^{\ast}=c_{H}^{-1},c_{H}^{\ast}=C_{H}^{-1}$ such that
\begin{equation}
c_{H}\mathbf{I}_{d_{2}J_{n}}\leq E( n^{-1}\mathbf{H}_{\mathbf{BB}})\leq
C_{H}\mathbf{I}%
_{d_{2}J_{n}}; \label{EQ:Bound_H&S}
\end{equation}
\textup{(ii)} with probability approaching $1$ as $n\rightarrow\infty$,
\begin{equation}
c_{H}\mathbf{I}_{d_{2}J_{n}}\leq n^{-1}\mathbf{H}_{\mathbf{BB}}\leq
C_{H}%
\mathbf{I}_{d_{2}J_{n} }. \label{EQ:Bound_Hhat}
\end{equation}
\end{lemma}%

Since the proof of Lemma~\ref{LEM:SandHposdef} is a little complicated,
we provide it in the supplemental article (Ma, Song and Wang~\cite{MSW11}). The
proofs of Lemmas~\ref{LEM:UhatU} to~\ref{LEM:varbetatilda} below are
also provided in (Ma, Song and Wang~\cite{MSW11}).

\begin{lemma}
\label{LEM:UhatU} Define $\widehat{\mathbf{U}}=( \sum_{i=1}^{n}%
\underline{\mathbf{B}}_{i}^{\T}\underline{\mathbf{X}}_{i})
_{d_{2}J_{n}\times d_{1} }$, where $\underline{\mathbf{B}}_{i}$ is
given in (\ref{EQ:B}). Under Assumptions \textup{(A1)--(A5)}, there exist constants $%
0<c_{U}<C_{U}<\infty$, such that with probability approaching $1$ as $%
n\rightarrow\infty$, $c_{U}\mathbf{I}_{d_{1}}\leq( n^{-1}h)
\widehat{\mathbf{U}}^{\T}\widehat{\mathbf{U}}\leq C_{U}\mathbf{I}_{d_{1}}.$
\end{lemma}

\begin{lemma}
\label{LEM:H11}Under Assumptions \textup{(A1)--(A5)}, there exist constants $%
0<c_{H_{1}}<C_{H_{1}}<\infty$, such that with probability approaching $1$
as $n\rightarrow\infty$, $c_{H_{1}}\mathbf{I}_{d_{1}}\leq n\mathbf{H}%
^{11}\leq C_{H_{1}}\mathbf{I}_{d_{1}}$, where $\mathbf{H}^{11}$ is
given in~\textup{(\ref{EQ:DVDinverse})}.
\end{lemma}

Let $\widetilde{\bolds{\beta}}_{\mu}$ and $\widetilde{\bolds{\beta}}_{e}$
be the solutions of (\ref{EQ:betahat}) with $\underline{\mathbf{Y}}_{i}$
replaced by $\underline{\bolds{\mu}}_{i}$ and $\underline{\mathbf{e}}%
_{i}=\underline{\mathbf{Y}}_{i}-\underline{\bolds{\mu}}_{i}$,
respectively. Then $\widehat{\bolds{\beta}}-\bolds{\beta}_{0}=(
\widetilde{\bolds{\beta}}_{\mu}-\bolds{\beta}_{0})+\widetilde{%
\bolds{\beta}}_{e}$.

\begin{lemma}
\label{LEM:betahatm-beta} Under Assumptions \textup{(A1)--(A5)}, $\Vert
\widetilde{%
\bolds{\beta}}_{\mu}-\bolds{\beta}_{0}\Vert=\mathrm{o}_{P}(
n^{-1/2}) $.
\end{lemma}

Note that $\widetilde{\bolds{\beta}}_{e}=\mathbf{H}^{11}\{
\sum_{i=1}^{n}\underline{\mathbf{X}}_{i}^{\T}\mathbf{V}_{i}^{-1}\mathbf{
\underline{e}}_{i}-\mathbf{H}_{\mathbf{XB}}\mathbf{H}_{\mathbf{BB}%
}^{-1}\sum_{i=1}^{n}\underline{\mathbf{B}}_{i}^{\T}\mathbf{V}_{i}^{-1}%
\underline{\mathbf{e}}_{i}\} $; thus we can show that the conditional
variance $\mathrm{%
\mathrm{Var}}( \widetilde{\bolds{\beta}}_{e}\vert\mathbb{X},%
\mathbb{Z}) $ equals
\begin{equation}
\mathbf{H}^{11}\sum_{i=1}^{n}\{ \underline{\mathbf{X}}_{i}-\underline{\mathbf{B}}_{i}\mathbf{H}_{\mathbf{BB}}^{-1}\mathbf{H}_{\mathbf{BX}%
}\} ^{\T}\mathbf{V}_{i}^{-1}\bolds{\Sigma}_{i}\mathbf{V}%
_{i}^{-1}\{ \underline{\mathbf{X}}_{i}-\underline{\mathbf{B}}_{i}%
\mathbf{H}_{\mathbf{BB}}^{-1}\mathbf{H}_{\mathbf{BX}}\} \mathbf{H}%
^{11}. \label{EQ:varbetahat}
\end{equation}

\begin{lemma}
\label{LEM:betahateasymp}Under Assumptions \textup{(A1)--(A5)}, as $n\rightarrow
\infty$,
\[
\{ \mathrm{Var}( \widetilde{\bolds{\beta}}_{e}\vert\mathbb{%
X},\mathbb{Z} ) \} ^{-1/2}( \widetilde{\bolds{\beta
}}_{e}) \longrightarrow N( 0,\mathbf{I}_{d_{1} }) .
\]
\end{lemma}

\begin{lemma}
\label{LEM:varbetatilda}Under Assumptions \textup{(A1)--(A5)}, for the covariance
matrix $\bolds{\Omega}( \mathbb{V},\bbSigma ) $ defined
in (\ref{EQ:IV}), $c_{V}^{\ast}\mathbf{I}_{d_{1}}\leq\bolds{\Omega}%
( \mathbb{V},\bbSigma ) \bolds{\leq}C_{V}^{\ast}\mathbf{I}_{d_{1}}$ and $\mathrm{Var}( \widetilde{\bolds{\beta}}%
_{e}\vert\mathbb{X},\mathbb{Z} ) =n^{-1}\bolds{\Omega}%
( \mathbb{V},\bbSigma ) +\mathrm{O}_{P}(
n^{-3/2}+n^{-1}h^{2p}) $.
\end{lemma}

Theorem~\ref{THM:betahat-normality} follows from Lemmas \ref%
{LEM:betahatm-beta},~\ref{LEM:betahateasymp} and~\ref{LEM:varbetatilda}.

\subsection{\texorpdfstring{Proof of Theorem \protect\ref{THM:asymptoticphihat}}{Proof of Theorem 2}}

From (\ref{EQ:betagammahat}) and (\ref{EQ:DVDinverse}), we obtain
\begin{equation}
\widehat{\bolds{\gamma}}=\mathbf{H}^{22}\Biggl( \sum_{i=1}^{n}\underline{\mathbf{B}}_{i}^{\T}\mathbf{V}_{i}^{-1}\underline
{\mathbf{Y}}_{i}-\mathbf{H}_{\mathbf{BX}}\mathbf{H}_{\mathbf{XX}}^{-1}\sum_{i=1}^{n}\underline{\mathbf{X}}_{i}^{\T}\mathbf{V}_{i}^{-1}\underline{\mathbf{Y}}_{i}\Biggr)
. \label{EQ:lamdahat}
\end{equation}
Following the same idea as that in the proof of Lemma~\ref{LEM:H11}, we
have that
there exist constants $0<c_{H_{2}}<C_{H_{2}}<\infty$, such that with
probability approaching $1$ as $n\rightarrow\infty$, $c_{H_{2}}\mathbf
{I}%
_{d_{2}J_{n}}\leq n\mathbf{H}^{22}\leq C_{H_{2}}\mathbf{I}_{d_{2}J_{n}}$.
Letting $\widetilde{\bolds{\gamma}}_{\mu}$ and $\widetilde{\bolds
{\gamma}%
}_{e}$ be the solutions of (\ref{EQ:lamdahat}) with $\underline{\mathbf
{Y}}%
_{i}$ replaced by $\underline{\bolds{\mu}}_{i}$ and $\mathbf
{\underline{{e}%
}}_{i}=\underline{\mathbf{Y}}_{i}-\underline{\bolds{\mu}}_{i}$,
respectively, $\widehat{\bolds{\gamma}}-\bolds{\gamma}=( \widetilde{
\bolds{\gamma}}_{\mu}-\bolds{\gamma}) +\widetilde{\bolds{\gamma}}%
_{e}$. Letting $\widehat{\Pi}_{n,\mathbf{X}}$ be the projection on $\{
\underline{\mathbf{X}}_{i}\} _{i=1}^{n}$ to the empirical inner
product, $\widetilde{\bolds{\gamma}}_{\mu}-\bolds{\gamma}$ equals
\begin{eqnarray*}
&& \mathbf{H}^{22}\Biggl[ \sum_{i=1}^{n}\underline{\mathbf{B}}_{i}^{\T}\mathbf
{%
V}_{i}^{-1}\Biggl\{ \sum_{l=1}^{d_{2}}\eta_{l}(\mathbf{Z}%
_{il})\Biggr\} -\mathbf{H}_{\mathbf{BX}}\mathbf{H}_{\mathbf{XX}%
}^{-1}\sum_{i=1}^{n}\underline{\mathbf{X}}_{i}^{\T}\mathbf{V}%
_{i}^{-1}\Biggl\{ \sum_{l=1}^{d_{2}}\eta_{l}(\mathbf{Z}_{il})\Biggr\} %
\Biggr] -\bolds{\gamma} \\
&&\quad =\mathbf{H}^{22}\sum_{i=1}^{n}\underline{\mathbf{B}}_{i}^{\T}\mathbf
{V}%
_{i}^{-1}\Biggl[ \Biggl\{ \sum_{l=1}^{d_{2}}\eta_{l}(\mathbf{Z}_{il})-%
\underline{\mathbf{B}}_{i}\bolds{\gamma}\Biggr\} -\widehat{\Pi}_{n,%
\mathbf{X}}\Biggl\{ \sum_{l=1}^{d_{2}}\eta_{l}(\mathbf{Z}_{il})-%
\underline{\mathbf{B}}_{i}\bolds{\gamma}\Biggr\} \Biggr] \\
&&\quad=n\mathbf{H}^{22}%
\mathbf{S},
\end{eqnarray*}
where $\mathbf{S}=( S_{11},\ldots,S_{J_{n}d_{2}}) $, with
\[
S_{s,l}=n^{-1}\sum_{i=1}^{n}\bigl( \mathbf{B}_{i}^{( s,l)
}\bigr) ^{\T}\mathbf{V}_{i}^{-1}\Biggl[ \Biggl\{
\sum_{l=1}^{d_{2}}\eta_{l}(\mathbf{Z}_{il})-\underline{\mathbf{B}}%
_{i}\bolds{\gamma}\Biggr\} -\widehat{\Pi}_{n,\mathbf{X}}\Biggl\{
\sum_{l=1}^{d_{2}}\eta_{l}(\mathbf{Z}_{il})-\underline{\mathbf{B}}%
_{i}\bolds{\gamma}\Biggr\} \Biggr] ,
\]
and $\mathbf{B}_{i}^{( s,l) }=[ \{ B_{s,l}(
Z_{i1l}) ,\ldots,B_{s,l}( Z_{im_{i}l}) \} ^{\T}]
_{m_{i}\times1}$. Let $\Delta\eta(\underline{\mathbf{Z}}_{i})=$ $%
\sum_{l=1}^{d_{2}}\eta_{l}(\mathbf{Z}_{il})-\underline{\mathbf{B}}%
_{i}\bolds{\gamma}$, then the Cauchy--Schwarz inequality implies that
\[
\vert S_{s,l}\vert\leq\Biggl\{ n^{-1}\sum_{i=1}^{n}\bigl(
\mathbf{B}_{i}^{( s,l) }\bigr) ^{\T}\mathbf{V}_{i}^{-1}\mathbf{B}%
_{i}^{( s,l) }\Biggr\} ^{1/2}\Vert\Delta\eta-\widehat{\Pi
}_{n,\mathbf{X}}( \Delta\eta) \Vert_{n}=\mathrm{O}_{P}(
h^{p}) ,
\]
thus $\Vert\widetilde{\bolds{\gamma}}_{\mu}-\bolds{\gamma}%
\Vert=$ $\mathrm{O}_{P}( J_{n}^{1/2}h^{p}) $. For any $\mathbf{c}%
\in\mathcal{R}^{J_{n}d_{2}}$ with $\Vert\mathbf{c}\Vert=1$, we
write $\mathbf{c}^{\T}\widetilde{\bolds{\gamma}}_{e}=\sum%
_{i=1}^{n}a_{i}\varepsilon_{i}$, \ where $\varepsilon_{i}$ are
independent conditioning on $( \mathbb{X},\mathbb{Z}) $ and
\[
a_{i}^{2}=\mathbf{c}^{\T}\mathbf{H}^{22}\{ \underline{\mathbf{B}}_{i}-%
\underline{\mathbf{X}}_{i}\mathbf{H}_{\mathbf{XX}}^{-1}\mathbf
{H}_{\mathbf{XB%
}}\} ^{\T}\mathbf{V}_{i}^{-1}\bolds{\Sigma}_{i}\mathbf{V}%
_{i}^{-1}\{ \underline{\mathbf{B}}_{i}-\underline{\mathbf{X}}_{i}%
\mathbf{H}_{\mathbf{XX}}^{-1}\mathbf{H}_{\mathbf{XB}}\} \mathbf{H}^{22}%
\mathbf{c}.
\]
Following the same arguments as those in Lemma~\ref{LEM:betahateasymp},
we have $\max_{1\leq i\leq n}\vert a_{i}\vert=\mathrm{O}_{P}(
J_{n}^{1/2}n^{-1}) $. Thus $\Vert\widetilde{\bolds{\gamma}}%
_{e}\Vert\leq J_{n}^{1/2}\vert\mathbf{c}^{\T}\widetilde{\bolds{%
\gamma}}_{e}\vert=J_{n}^{1/2}\vert
\sum_{i=1}^{n}a_{i}\varepsilon_{i}\vert=\mathrm{O}_{P}(
J_{n}^{1/2}n^{-1/2}) $. Therefore, $\Vert\widehat{\bolds{\gamma
}}_{l}-\bolds{\gamma}_{l}\Vert
= \mathrm{O}_{P}( J_{n}^{1/2}h^{p}+J_{n}^{1/2}n^{-1/2}) $. Because $\widehat{%
\eta}_{l}(z_{l})=\mathbf{B}_{l}^{\ast}( z_{l}) ^{\T}\widehat{%
\bolds{\gamma}}_{l}$, $\widetilde{\eta}_{l}(z_{l})=\mathbf
{B}_{l}^{\ast
}( z_{l}) ^{\T}\bolds{\gamma}_{l}$ and $\vert\widehat{\eta
}_{l}-\widetilde{\eta}_{l}\vert_{L_{2}}^{2}=\Vert\widehat{%
\bolds{\gamma}}_{l}-\bolds{\gamma}_{l}\Vert^{2}\times \mathrm{O}_{P}(
1) =\mathrm{O}_{P}( J_{n}h^{2p}+J_{n}n^{-1}) $. Thus one has
\[
\vert\widehat{\eta}_{l}-\eta_{l}\vert_{L_{2}}^{2}\leq
2( \vert\widehat{\eta}_{l}-\widetilde{\eta}_{l}\vert
_{L_{2}}^{2}+\vert\widetilde{\eta}_{l}-\eta_{l}\vert
_{L_{2}}^{2}) =\mathrm{O}_{P}( J_{n}h^{2p}+J_{n}n^{-1}).
\]

\subsection{\texorpdfstring{Proof of Theorem \protect\ref{THM:rootn}}{Proof of Theorem 3}}

Let $\tau_{n}=n^{-1/2}+a_{n}$. It suffices to show that for any given $
\zeta>0$, there exists a large constant $C$ such that
\begin{equation}
P\Bigl\{\sup_{\Vert\mathbf{u}\Vert=C}Q_{\mathcal{P}}(\bolds{\beta}%
_{0}+\tau_{n}\mathbf{u})>Q_{\mathcal{P}}(\bolds{\beta}_{0})\Bigr\}\geq
1-\zeta. \label{EQ:existence}
\end{equation}
Plugging $\bolds{\gamma}( \bolds{\beta}) $ in (\ref%
{EQ:gamma(beta)}) into $Q( \bolds{\beta}) $ defined in (\ref%
{DEF:Qbeta}), we have
\begin{eqnarray*}
Q( \bolds{\beta}) &=&\frac{1}{2}\sum_{i=1}^{n}\Biggl[ \mathbf{Y}%
_{i}-\Biggl\{ \underline{\mathbf{X}}_{i}\bolds{\beta}+\underline{\mathbf{B}}
_{i}\mathbf{H}_{\mathbf{BB}}^{-1}\sum_{i=1}^{n}\mathbf{B}_{i}^{\T
}\mathbf{V}%
_{i}^{-1}( \mathbf{Y}_{i}-\underline{\mathbf{X}}_{i}\bolds{\beta}%
) \Biggr\} \Biggr] ^{\T} \\
&&\phantom{\frac{1}{2}\sum_{i=1}^{n}}{}\times\mathbf{V}_{i}^{-1}\Biggl[ \mathbf{Y}_{i}-\Biggl\{ \underline{\mathbf{X}}_{i}\bolds{\beta}+\underline{\mathbf{B}}_{i}\mathbf{H}_{\mathbf{BB}
}^{-1}\sum_{i=1}^{n}\underline{\mathbf{B}}_{i}^{\T}\mathbf{V}_{i}^{-1}(
\mathbf{Y}_{i}-\underline{\mathbf{X}}_{i}\bolds{\beta}) \Biggr\} %
\Biggr] .
\end{eqnarray*}
Thus $Q( \bolds{\beta}) =\frac{1}{2}\sum_{i=1}^{n}(
\mathbf{Y}_{i}-\widehat{\underline{\mathbf{X}}}_{i}\bolds{\beta
}-\widehat{%
\Pi}_{n}\mathbf{Y}_{i}) ^{\T}\mathbf{V}_{i}^{-1}( \mathbf{Y}_{i}-%
\widehat{\underline{\mathbf{X}}}_{i}\bolds{\beta}-\widehat{\Pi}_{n}%
\mathbf{Y}_{i}) . $ Let $U_{n,1}=Q(\bolds{\beta}_{0}+\tau_{n}%
\mathbf{u})-Q(\bolds{\beta}_{0})$ and $U_{n,2}=n_{\T}\sum_{k=1}^{r}\{
p_{%
\lambda_{k}}(|\beta_{k0}+\tau_{n}u_{k}|)-p_{\lambda_{k}}(|\beta
_{k0}|)\} $, where $r$ is the number of components of $\bolds{\beta
}_{10}$%
. Note that $p_{\lambda_{k}}( 0) =0$ and $p_{\lambda_{k}}(
|\beta|) \geq0$ for all $\beta$. Thus, $Q_{\mathcal{P}}(\bolds{%
\beta}_{0}+\tau_{n}\mathbf{u})-Q_{\mathcal{P}}(\bolds{\beta
}_{0})\geq
U_{n,1}+U_{n,2}$.

For $U_{n,1}$, we have $Q(\bolds{\beta}_{0}+\tau_{n}\mathbf
{u})=Q(\bolds{%
\beta}_{0})+\tau_{n}\mathbf{u}^{\T}\dot{Q}(\bolds{\beta}_{0})+\frac
{1}{2}%
\tau_{n}^{2}\mathbf{u}^{\T}\ddot{Q}(\bolds{\beta}^{\ast})\mathbf
{u}, $
where $\ddot{Q}(\bolds{\beta})=\sum_{i=1}^{n}\widehat{\underline
{\mathbf{X}%
}}_{i}^{\T}\mathbf{V}_{i}^{-1}\widehat{\underline{\mathbf{X}}}_{i}$, $%
\bolds{\beta}^{\ast}=t( \bolds{\beta}_{0}+n^{-1/2}\mathbf{u}) +(
1-t) \bolds{\beta}_{0}$, $t\in\lbrack0,1]$. Note that
\begin{eqnarray*}
\dot{Q}(\bolds{\beta}_{0}) &=&\sum_{i=1}^{n}\widehat{\underline
{\mathbf{X}}%
}_{i}^{\T}\mathbf{V}_{i}^{-1}( \mathbf{Y}_{i}-\widehat{\underline{%
\mathbf{X}}}_{i}\bolds{\beta}_{0}\mathbf{-}\widehat{\Pi}_{n}\mathbf
{Y}_{i}%
) \\
&=&\sum_{i=1}^{n}\widehat{\underline{\mathbf{X}}}_{i}^{\T}\mathbf{V}%
_{i}^{-1}\Biggl\{ \sum_{l=1}^{d_{2}}\eta_{l}( \mathbf{Z}_{il}) -%
\widehat{\Pi}_{n}\sum_{l=1}^{d_{2}}\eta_{l}( \mathbf{Z}_{il})
\Biggr\} +\sum_{i=1}^{n}\widehat{\underline{\mathbf{X}}}_{i}^{\T}\mathbf{V}%
_{i}^{-1}( \underline{\mathbf{e}}_{i}-\widehat{\Pi}_{n}\underline{\mathbf{e}}_{i}) ,
\end{eqnarray*}
where $\underline{\mathbf{e}}_{i}=\underline{\mathbf{Y}}_{i}-\underline{\bolds{\mu}}_{i}$. Mimicking the proof for Lemmas \ref%
{LEM:betahatm-beta} and~\ref{LEM:betahateasymp}, we have
\begin{eqnarray*}
\sum_{i=1}^{n}\widehat{\underline{\mathbf{X}}}_{i}^{\T}\mathbf{V}%
_{i}^{-1}\Biggl\{ \sum_{l=1}^{d_{2}}\eta_{l}( \mathbf{Z}_{il}) -%
\widehat{\Pi}_{n}\sum_{l=1}^{d_{2}}\eta_{l}( \mathbf{Z}_{il})
\Biggr\} &=&\mathrm{o}_{P}( n^{1/2}) ,
\\
\sum_{i=1}^{n}\widehat{\underline{\mathbf{X}}}_{i}^{\T}\mathbf{V}%
_{i}^{-1}( \underline{\mathbf{e}}_{i}-\widehat{\Pi}_{n}\underline{\mathbf{e}}_{i}) &=&\mathrm{O}_{P}(n^{1/2}) .
\end{eqnarray*}
Thus $\tau_{n}\mathbf{u}^{\T}\dot{Q}(\bolds{\beta
}_{0})=\mathrm{O}_{P}(n^{1/2}\tau
_{n})\Vert\mathbf{u}\Vert$. By the proof of Lemma~\ref{LEM:H11}%
, we obtain that $\frac{1}{2}\tau_{n}^{2}\mathbf{u}^{\T}\times \ddot
{Q}(\bolds{%
\beta}_{0})\mathbf{u}=\mathrm{O}_{P}(n\tau_{n}^{2})+\mathrm{o}_{P}(1)$. Thus
\begin{equation}
U_{n,1}=\mathrm{O}_{P}(n^{1/2}\tau_{n})+\mathrm{O}_{P}(n\tau_{n}^{2})+\mathrm{o}_{P}(1).
\label{EQ:Dn1}
\end{equation}
For $U_{n,2}$, by a Taylor expansion,
\[
p_{\lambda_{k}}(|\beta_{k0}+\tau_{n}u_{k}|)=p_{\lambda_{k}}(|\beta
_{k0}|)+\tau_{n}u_{k}p_{\lambda_{k}}^{\prime}( \vert\beta
_{k0}\vert) \mathrm{sgn}( \beta_{k0}) +\tfrac{1}{2}%
\tau_{n}^{2}u_{k}^{2}p_{\lambda_{k}}^{\prime\prime}( \vert
\beta_{k}^{\ast}\vert) ,
\]
where $\beta_{k}^{\ast}=(1-t)\beta_{k0}+t(\beta
_{k0}+n^{-1/2}u_{k})$, $%
t\in[0,1]$ and
\[
p_{\lambda_{k}}(|\beta_{k0}+\tau_{n}u_{k}|)=p_{\lambda_{k}}(|\beta
_{k0}|)+\tau_{n}u_{k}p_{\lambda_{k}}^{\prime}( \vert\beta_{k0}\vert
) \mathrm{sgn}( \beta_{k0}) +\tfrac{1}{2}\tau
_{n}^{2}u_{k}^{2}p_{\lambda_{k}}^{\prime
\prime}( \vert\beta_{k0}\vert)+\mathrm{o}( n^{-1}).
\]
Thus, by the Cauchy--Schwarz inequality,
\begin{eqnarray*}
n_{\T}^{-1}U_{n,2} &=&\tau_{n}\sum_{k=1}^{r}u_{k}p_{\lambda
_{k}}^{\prime
}( \vert\beta_{k0}\vert) \mathrm{sgn}( \beta
_{k0}) +\frac{1}{2}\tau_{n}^{2}\sum_{k=1}^{r}u_{k}^{2}p_{\lambda
_{k}}^{\prime\prime}( \vert\beta_{k0}\vert) \\
&\leq&\sqrt{r}\tau_{n}a_{n}\Vert\mathbf{u}\Vert+\frac{1}{2}\tau
_{n}^{2}w_{n}\Vert\mathbf{u}\Vert^{2}=C\tau_{n}^{2}\bigl(\sqrt{r}+w_{n}C\bigr).
\end{eqnarray*}
As $w_{n}\rightarrow0$, the first two terms on the right-hand side of
(\ref%
{EQ:Dn1}) dominate $U_{n,2}$ by taking $C$ sufficiently large. Hence
(\ref%
{EQ:existence}) holds for sufficiently large $C$.

\subsection{\texorpdfstring{Proof of Theorem \protect\ref{THM:oracle}}{Proof of Theorem 4}}

We first show that the estimator $\widehat{\bolds{\beta}}^{\P}$ must
possess the sparsity property $\widehat{\bolds{\beta}}_{2}=0$, which is
stated as follows.

\begin{lemma}
\label{LEM:sparsity} Under the conditions of Theorem~\ref{THM:oracle}, with
probability tending to 1, for any given $\bolds{\beta}_{1}$ satisfying
that $\Vert\bolds{\beta}_{1}-\bolds{\beta}_{10}\Vert=\mathrm{O}_{P}(n^{-1/2})$
and any constant $C$,
\[
Q_{\mathcal{P}}\{(\bolds{\beta}_{1}^{\T},\mathbf{0}^{\T})^{\T}\}=\min
_{\Vert\bolds{%
\beta}_{2}\Vert\leq Cn^{-1/2}}Q_{\mathcal{P}}\{(\bolds{\beta
}_{1}^{\T},%
\bolds{\beta}_{2}^{\T})\}.
\]
\end{lemma}

\begin{pf}
To prove that the maximizer is obtained at $\bolds{\beta}_{2}=0$, it
suffices to show that with probability tending to 1, as $n\rightarrow
\infty
$, for any $\bolds{\beta}_{1}$ satisfying $\Vert\bolds{\beta}_{1}-%
\bolds{\beta}_{10}\Vert=\mathrm{O}_{P}(n^{-1/2})$, and $\Vert\bolds{\beta}%
_{2}\Vert\leq Cn^{-1/2}$, $\partial Q_{\mathcal{P}}(\bolds{\beta}%
)/\partial\beta_{k}$ and $\beta_{k}$ have different signs for $\beta
_{k}\in(-Cn^{-1/2},Cn^{-1/2})$, for $k=r+1,\ldots,d_{1}$. Note that
\[
\dot{Q}_{\mathcal{P},k}( \bolds{\beta}) \equiv\frac{\partial
Q_{\mathcal{P}}(\bolds{\beta})}{\partial\beta_{k}}=\dot{Q}_{k}(
\bolds{\beta}) +n_{\T}p_{\lambda_{kn}}^{\prime}( \vert
\beta_{k}\vert) \sgn(\beta_{k}),
\]
where $\dot{Q}_{k}( \bolds{\beta}) =\dot{Q}_{k}( \bolds{%
\beta}_{0}) +\sum_{k^{\prime}=1}^{d_1}\ddot{Q}_{kk^{\prime}}\{t%
\beta_{k^{\prime}}+(1-t)\beta_{0k^{\prime}}\}( \beta
_{k^{\prime}}-\beta_{0k^{\prime}}) $, $t\in[0,1]$,
\[
\dot{Q}_{k}( \bolds{\beta}_{0}) =e_{k}^{\T}\sum_{i=1}^{n}%
\widehat{\underline{\mathbf{X}}}_{i}^{\T}\mathbf{V}_{i}^{-1}( \mathbf{Y}
_{i}-\widehat{\underline{\mathbf{X}}}_{i}\bolds{\beta}_{0}-\widehat
{\Pi}%
_{n}\mathbf{Y}_{i}) .
\]
It follows by the similar arguments as given in the proofs of Theorems
\ref%
{THM:betahat-normality} and~\ref{THM:rootn} that
\begin{eqnarray*}
\dot{Q}_{k}( \bolds{\beta}_{0}) &=&e_{k}^{\T}\sum_{i=1}^{n}%
\widehat{\underline{\mathbf{X}}}_{i}^{\T}\mathbf{V}_{i}^{-1}\Biggl\{
\sum_{l=1}^{d_{2}}\eta_{l}( \mathbf{Z}_{il}) -\widehat{\Pi}%
_{n}\sum_{l=1}^{d_{2}}\eta_{l}( \mathbf{Z}_{il}) \Biggr\} +e_{k}^{\T}\sum
_{i=1}^{n}\widehat{\underline{\mathbf{X}}}_{i}^{\T}\mathbf{V%
}_{i}^{-1}( \underline{\mathbf{e}}_{i}-\widehat{\Pi}_{n}\underline{\mathbf{e}}_{i}) \\
&=&n\Biggl\{ n^{-1}\sum_{i=1}^{n}\Xi_{k}(\underline{\mathbf{Y}}_{i},%
\underline{\mathbf{X}}_{i},\bolds{\underline
{Z}}_{i})+\mathrm{o}_{P}(n^{-1/2})\Biggr\} ,
\end{eqnarray*}
where $\Xi_{k}(\underline{\mathbf{Y}}_{i},\underline{\mathbf{X}}_{i},%
\underline{\mathbf{Z}}_{i})$ is the $k$th element of matrix $\widehat{%
\underline{\mathbf{X}}}_{i}^{\T}\mathbf{V}_{i}^{-1}( \underline{\mathbf{
e}}_{i}-\widehat{\Pi}_{n}\underline{\mathbf{e}}_{i}) $. According to
Lemma~\ref{LEM:varbetatilda}, we have
\begin{eqnarray*}
n^{-1}\ddot{Q}( \bolds{\beta}_{0})&=&E\Biggl( n^{-1}\sum_{i=1}^{n}\widetilde
{\underline{%
\mathbf{X}}}_{i}^{\T}\mathbf{V}_{i}^{-1}\widetilde{\underline{\mathbf
{X}}}%
_{i}\Biggr)+\mathrm{o}_{P}( 1) =\mathbf{R}+\mathrm{o}_{P}( 1),
\\
\frac{1}{n}\sum_{k^{\prime}=1}^{d_1}\ddot{Q}_{kk^{\prime}}( \beta
_{k^{\prime}}-\beta_{0k^{\prime}}) &=&(\bolds{\beta}-\bolds{\beta}%
_{0})^{\T}\bigl( R_{k}+\mathrm{o}_{P}( 1) \bigr) ,
\end{eqnarray*}
where $R_{k}$ is the $k$th column of $\mathbf{R}$. Note that $\Vert
\bolds{%
\beta}-\bolds{\beta}_{0}\Vert=\mathrm{O}_{P}(n^{-1/2})$ by the assumption.
Thus, $%
n^{-1}\dot{Q}_{k}( \bolds{\beta}) $ is of the order $%
\mathrm{O}_{P}(n^{-1/2})$. Therefore, for any nonzero $\beta_{k}$ and
$k=r+1,\ldots
,d_{1}$,
\[
\dot{Q}_{\mathcal{P},k}( \bolds{\beta})=n\lambda_{kn}\biggl\{\lambda
_{kn}^{-1}p_{\lambda_{kn}}^{\prime}(
\vert\beta_{k}\vert) \sgn(\beta_{k})+\mathrm{O}_{P} \biggl(\frac{1%
}{\sqrt{n}\lambda_{kn}} \biggr)\biggr\}.
\]
Since $\liminf_{n\rightarrow\infty}\liminf_{\beta_{k}\rightarrow
0^{+}}\lambda_{kn}^{-1}p_{\lambda_{kn}}^{\prime}(|\beta_{k}|)>0$
and $%
\sqrt{n}\lambda_{kn}\rightarrow\infty$, the sign of the derivative is
determined by that of $\beta_{k}$. Thus the desired result is obtained.
\end{pf}

\begin{pf*}{Proof of Theorem \protect\ref{THM:oracle}}
From Lemma~\ref{LEM:sparsity}, it follows that $\widehat{\bolds{\beta
}}%
_{2}^{\P}=\mathbf{0}$.
\begin{eqnarray*}
\dot{Q}_{\mathcal{P}}( \bolds{\beta}) &=&\dot{Q}(\bolds{\beta}%
_{0})+\ddot{Q}(\bolds{\beta}^{\ast})( \bolds{\beta}-\bolds{\beta}%
_{0}) +n_{\T}\{ p_{\lambda_{kn}}^{\prime}( \vert\beta
_{k0}\vert) \operatorname{sign}( \beta_{k0}) \}
_{k=1}^{r} \\
&&{}+\Biggl\{ \sum_{k=1}^{r}p_{\lambda_{kn}}^{\prime\prime}( \vert
\beta_{k0}\vert) +\mathrm{o}_{P}( 1) \Biggr\} ( \widehat{\bolds{\beta}}\,_{k1}^{\P}-\beta_{k0}) ,
\end{eqnarray*}
where $\bolds{\beta}^{\ast}=t\bolds{\beta}_{0}+( 1-t)
\bolds{\beta}$, $t\in\lbrack0,1]$. Using an argument similar to the
proof of Theorem~\ref{THM:rootn}, it can be shown that there exists a $%
\widehat{\bolds{\beta}}_{1}^{\P}$ in Theorem~\ref{THM:rootn} that is
a root-$n$
consistent local minimizer of $Q_{\mathcal{P}}\{( \bolds{\beta}_{1}
^{\T},%
\mathbf{0}^{^{\T}}) ^{\T}\} $%
, satisfying the penalized least squares equations $\dot{Q}_{\mathcal
{P}}%
[\{ ( \widehat{\bolds{\beta}}_{1}^{\P}) ^{\T},%
\mathbf{0}^{^{\T}}\} ^{\T}] =\mathbf{0}$. Mimicking the proofs
for Lemmas~\ref{LEM:betahatm-beta} and~\ref{LEM:betahateasymp} indicates
that the left hand side of the above equation can be written as
\begin{eqnarray*}
&&\hspace*{-4pt}
n^{-1}\sum_{i=1}^{n}\widehat{\underline{\mathbf{X}}}_{1i}^{\T}%
\mathbf{V}_{i}^{-1}( \underline{\mathbf{e}}_{i}-\hat{\Pi}_{n}\underline{\mathbf{e}}_{i}) +\{ p_{\lambda_{kn}}^{\prime}(
\vert\beta_{k0}\vert) \operatorname{sign}( \beta
_{k0}) \} _{k=1}^{r}+\mathrm{o}_{P}( n^{-1/2})
\\
&&\hspace*{-5pt}\quad{}+\Biggl\{ E\Biggl( n^{-1}\sum_{i=1}^{n}%
\widetilde{\underline{\mathbf{X}}}_{1i}^{\T}\mathbf
{V}_{i}^{-1}\widetilde{%
\underline{\mathbf{X}}}_{1i}\Biggr) +\mathrm{o}_{P}( 1) \Biggr\}(
\widehat{\bolds{\beta}}_{1}^{\P}-\bolds{\beta}_{10})\\
&&\hspace*{-5pt}\quad{}+\Biggl\{ \sum
_{k=1}^{r}p_{\lambda_{kn}}^{\prime\prime}( \vert
\beta_{k0}\vert) +\mathrm{o}_{P}( 1) \Biggr\}(
\widehat{\bolds{\beta}}_{1}^{\P}-\bolds{\beta}_{10}) .
\end{eqnarray*}
Thus we have
\begin{eqnarray*}
\mathbf{0} &=&n^{-1}\sum_{i=1}^{n}\widehat{\bolds{\underline
{X}}}_{1i}^{\T}%
\mathbf{V}_{i}^{-1}( \underline{\mathbf{e}}_{i}-\hat{\Pi}_{n}\underline{\mathbf{e}}_{i}) +\kappa_{n}+\mathrm{o}_{P}(n^{-1/2}) \\
&&{}+\Biggl\{ E\Biggl( n^{-1}\sum_{i=1}^{n}\widetilde{\underline{\mathbf{X}}}%
_{1i}^{\T}\mathbf{V}_{i}^{-1}\widetilde{\underline{\mathbf{X}}}_{1i}\Biggr) +%
\bolds{\Sigma}_{\lambda}+\mathrm{o}_{P}(1) \Biggr\} (\widehat{\bolds{%
\beta}}_{1}^{\P}-\bolds{\beta}_{10}).
\end{eqnarray*}
Similar arguments to Lemmas~\ref{LEM:betahateasymp} and \ref%
{LEM:varbetatilda} yield the asymptotic normality.
\end{pf*}

\end{appendix}

\section*{Acknowledgments}
Ma's research was supported by a dissertation fellowship from Michigan
State University. Wang's research was supported in part by NSF award
DMS-0905730.
The authors are grateful for the insightful comments from the editor,
an associate editor and anonymous referees.

\begin{supplement} [id=suppA]
\stitle{Supplement to ``Simultaneous variable selection and estimation
in semiparametric modeling of longitudinal/clustered data''}
\slink[doi]{10.3150/11-BEJ386SUPP}
\sdatatype{.pdf}
\sfilename{bej386\_supp.pdf}
\sdescription{We provide detailed proofs of Lemmas \ref
{LEM:SandHposdef} to~\ref{LEM:varbetatilda} stated in the \hyperref[appendix]{Appendix}.}
\end{supplement}

%

%

\printhistory

\end{document}